\documentclass[12pt,oneside,brazil,english]{amsart}
\usepackage[T1]{fontenc}
\usepackage[latin9]{inputenc}
\usepackage[letterpaper]{geometry}
\usepackage{float}
\usepackage{amsthm}
\usepackage{graphicx}
\usepackage{setspace}
\usepackage{amssymb}
\usepackage{esint}
\makeatletter
\numberwithin{equation}{section}
\numberwithin{figure}{section}
\theoremstyle{plain}
\newtheorem{thm}{Theorem}
  \theoremstyle{definition}
  \newtheorem{defn}[thm]{Definition}
  \theoremstyle{plain}
  \newtheorem{prop}[thm]{Proposition}
  \theoremstyle{plain}
  \newtheorem{cor}[thm]{Corollary}
  \theoremstyle{plain}
  \newtheorem{lem}[thm]{Lemma}

\makeatother

\usepackage{babel}

\begin{document}

\title{$K$-Theory of Boutet de Monvel algebras with classical SG-symbols on
the half space.}

\author{Pedro T. P. Lopes}
\address{Departamento de Matem\'atica da Universidade Federal de S\~ao Carlos,  13565-905, S\~ao Carlos, SP, Brazil}
\email{dritao@yahoo.com; pedrolopes@dm.ufscar.br}
\thanks{The authors were supported by CNPq, Brazil.}

\author{Severino T. Melo}
\address{Instituto de Matem\'atica e Estat\'istica, Universidade de S\~ao Paulo,  Rua do Mat\~ao 1010, 05508-090, S\~ao Paulo, SP, Brazil}
\email{toscano@ime.usp.br}

%    General info

\date{\today}

\begin{abstract}
We compute the $K$-groups of the $C^{*}$-algebra of bounded operators
generated by the Boutet de Monvel operators with classical SG-symbols of
order (0,0) and type 0 on $\mathbb{R}_{+}^{n}$, as defined by Schrohe,
Kapanadze and Schulze. In order to adapt the techniques used in Melo, Nest, Schick and Schrohe's work on the K-theory of Boutet de Monvel's algebra on compact manifolds, we regard the symbols as functions defined on the radial compactifications of $\mathbb{R}_{+}^{n}\times\mathbb{R}^{n}$ and $\mathbb{R}^{n-1}\times\mathbb{R}^{n-1}$.
This allows us to give useful descriptions of the kernel and the image of the continuous extension of the boundary
principal symbol map, which defines a $C^{*}$-algebra homomorphism. We are then able to compute the $K$-groups of the algebra using the standard K-theory six-term cyclic exact sequence associated to that homomorphism.

\end{abstract}
\maketitle
\tableofcontents{}

Boutet de Monvel introduced his algebra of pseudodifferential boundary value problems in 1971 \cite{boutetMonvel}.
It is essentially the smallest algebra that contains the differential elliptic
boundary value problems as well as their parametrices. In \cite{boutetMonvel}, using the Atiyah-Singer index Theorem \cite{AS1}, he proved that, in the case of compact manifolds with boundary, the Fredholm index for the elliptic elements of his algebra also factors through the K-theory of the cotangent bundle of the interior. Later that fact was used by several authors \cite{F,RempelSchulze} to give formulas for the index in terms of the symbols. 

In recent years, $C^{*}$-algebra $K$-theory techniques allowed
the extension and simplification of some results obtained by Boutet de Monvel \cite{MeloSchrohe,gaarde,ToscanoSchick,families,survey}.
In the context of noncompact manifolds, growth estimates naturally have to be imposed on the symbols. Using what became known as SG-estimates (whose definition goes back at least to Cordes \cite{parametrix} and Parenti \cite{Parenti}), Schrohe was able to construct a Boutet de Monvel
algebra on non-compact manifolds and to obtain results about the spectral
invariance of those operators \cite{SchroheSGboutet}. The SG-Boutet de Monvel algebra
for classical symbols (i.e., with symbols possessing asymptotic expansions in homogeneous components) was studied more recently by Kapanadze and Schulze
with operator-valued symbols techniques \cite{KapanadzeSchulze,KapanadzeSchulzearticle}.

Our study of the K-theory of the classical SG-Boutet de Monvel algebra on the half-space combines techniques used in the study of the K-theory of Boutet de Monvel's algebra on compact manifolds \cite{MeloSchrohe,ToscanoSchick} with the geometric characterization of SG-symbols via the radial compactification, as explained by Melrose \cite{Melrosescattering}, and used by Nicola \cite{Nicolakteoria} to compute the K-theory of the SG-algebra on $\mathbb{R}^n$.
Our strategy is  to find explicit descriptions of the image and the kernel of the $C^*$-homomorphism defined by the continuous extension of the (operator-valued) boundary principal symbol. As in \cite{Ditsche,MeloSchrohe,Cintia}, it then turns out that the analysis of the standard K-theory six-term cyclic exact sequence associated to that homomorphism suffices for the computation of the K-groups.

Many authors regard the leading terms in the two asymptotic expansions (with respect to $x$ and with respect to $\xi$) as two different principal symbols. Not only for $K$-theory computations, but for the study of the Fredholm property as well, it is very convenient, and it is the point of view adopted here, to regard those two symbols as only one function defined on the infinity points of a certain compactification of the cotangent bundle of the underlying manifold, which is the half-space in our case. Similarly, we also regard the boundary principal symbol as only one (operator-valued) function defined on the infinity points of a compactification of $\mathbb{R}^{n-1}\times\mathbb{R}^{n-1}$. This function is obtained by pointwise conjugation of the usual boundary principal symbols by a unitary operator valued function.

The proofs of many of the estimates needed in this paper are straightforward adaptations of published results, and would be considered standard by experts. We omit them, quoting works where similar proofs can be found. A very detailed exposition of our results can also be found in the first author's PhD thesis \cite{Tese}.

\section{Basic definitions.}
\label{sec:sec1}
In this section, we recall some of the main properties of the
SG-calculus. For proofs we refer to \cite{Schulze,Schrohecomplex,Rodinonicola,Parenti,KapanadzeSchulze,Cordes}.

We denote by $\left\langle \right\rangle :\mathbb{R}^{n}\to\mathbb{R}$ the function $\left\langle x\right\rangle :=\sqrt{1+\left|x\right|^{2}}$. The set $\{0,1,2,...\}$ of nonnegative integers is denoted by $\mathbb{N}_{0}$.
If $\mathcal{H}$ and $\tilde{\mathcal{H}}$ are Hilbert spaces,  then $\mathcal{B}(\mathcal{H},\tilde{\mathcal{H}})$ denotes the set of bounded operators from ${\mathcal{H}}$ to $\tilde{\mathcal{H}}$, $\mathcal{B}(\mathcal{H})$
denotes the set of bounded operators on $\mathcal{H}$ and $\mathcal{K}(\mathcal{H})$
denotes the set of compact operators on $\mathcal{H}$. As usual, $\mathcal{S}(\mathbb{R}^{n})$
denotes the Schwartz space of smooth functions whose derivatives are rapidly decreasing. 
Let $\Omega$ be an open set of $\mathbb{R}^{n}$, then $\mathcal{S}(\Omega)$ denotes the set of restrictions of Schwartz functions to $\Omega$. The main examples are $\mathcal{S}(\mathbb{R}_{+}^{n})$ and $\mathcal{S}(\mathbb{R}_{+}^{n}\times\mathbb{R}_{+}^{n})$, where 
$\mathbb{R}_{+}^{n}$ is the set $\{x\in\mathbb{R}^{n},\, x_{n}>0\}$.

Let us now recall the radial compactification, as introduced by Melrose \cite[section 6.3]{Melrosescattering}, see also \cite{Melrosenotas}. Let $\mathbb{S}_{+}^{n}=\{z\in\mathbb{R}^{n+1};\,|z|=1\,\mbox{and}\, z_{n+1}\ge0\}$.
The radial compactification 
is obtained regarding the map $RC:\mathbb{R}^{n}\to\mathbb{S}_{+}^{n}$, \[
RC(z)=\left(\frac{z}{\left\langle z\right\rangle },\frac{1}{\left\langle z\right\rangle }\right),\]
as an embedding.

Using this map we define:
\begin{defn}
\emph{(Classical SG-symbols)} The space of classical SG-symbols of order $(0,0)$,
denoted by $S_{cl}^{0,0}(\mathbb{R}^{n}\times\mathbb{R}^{n})$, is
the space of functions $a\in C^{\infty}(\mathbb{R}^{n}\times\mathbb{R}^{n})$
such that \[
a\circ\left(RC^{-1}\times RC^{-1}\right):\left(\mathbb{S}_{+}^{n}\cap\mathbb{R}_{+}^{n+1}\right)\times\left(\mathbb{S}_{+}^{n}\cap\mathbb{R}_{+}^{n+1}\right)\to\mathbb{C}\]
can be extended, uniquely, to a function $\tilde{a}\in C^{\infty}(\mathbb{S}_{+}^{n}\times\mathbb{S}_{+}^{n})$ - we note that    $\mathbb{S}_{+}^{n}\cap\mathbb{R}_{+}^{n+1}$ is the interior of the manifold with boundary $\mathbb{S}_{+}^{n}$.
Let $(\mu,\nu)\in\mathbb{R}^{2}$, then a classical SG-symbol of order $(\mu,\nu)\in\mathbb{R}^{2}$
is a function $a\in C^{\infty}(\mathbb{R}^{n}\times\mathbb{R}^{n})$
such that \[
(x,\xi)\mapsto\left\langle x\right\rangle ^{-\nu}\left\langle \xi\right\rangle ^{-\mu}a(x,\xi)\in S_{cl}^{0,0}(\mathbb{R}^{n}\times\mathbb{R}^{n}).\]
The set consisting of classical symbols of order $(\mu,\nu)$ is denoted by $S_{cl}^{\mu,\nu}(\mathbb{R}^{n}\times\mathbb{R}^{n})$

\end{defn}

It is necessary to make two remarks about this definition. The first one is that, as a consequence of the above definition, every function $a\in S_{cl}^{\mu,\nu}(\mathbb{R}^{n}\times\mathbb{R}^{n})$ satisfy the following estimates  \[\left|\partial_{x}^{\beta}\partial_{\xi}^{\alpha}a(x,\xi)\right|\le C_{\alpha\beta}\left\langle x\right\rangle ^{\nu-\left|\beta\right|}\left\langle \xi\right\rangle ^{\mu-\left|\alpha\right|},\]
where $C_{\alpha\beta}$ are positive constants that depend only on $a$, $\alpha$ and $\beta$. A smooth function that satisfies these estimates is called a SG-symbol. The set of SG-symbols is denoted by $S^{\mu,\nu}(\mathbb{R}^{n}\times\mathbb{R}^{n})$. It is clear that  $S_{cl}^{\mu,\nu}(\mathbb{R}^{n}\times\mathbb{R}^{n})\subset S^{\mu,\nu}(\mathbb{R}^{n}\times\mathbb{R}^{n})$.

As a second remark, we recall that a function $a$ is smooth on the manifold with corners $\mathbb{S}_{+}^{n}\times\mathbb{S}_{+}^{n}$ if, and only if, $a$ can be extended to a smooth function in $\mathbb{S}^{n}\times\mathbb{S}^{n}$. Using the variables $t=\frac{1}{|x|}$, $\Omega_{1}=\frac{x}{|x|}$ and $s=\frac{1}{|\xi|}$, $\Omega_{2}=\frac{\xi}{|\xi|}$, we see that $a\in S_{cl}^{0,0}(\mathbb{R}^{n}\times\mathbb{R}^{n})$ iff it has a smooth extension for $t=0$ and $s=0$, when it is written in terms of the variables $(x,s,\Omega_{2})$, $(t,\Omega_{1},\xi)$ and $(t,\Omega_{1},s,\Omega_{2})$. We can then take Taylor series in the variables $t$ and $s$ and conclude that the following asymptotic expansions for symbols of order $(0,0)$ hold:

\vspace{0,3cm}

(i) $a\sim\sum_{k=0}^{\infty}a_{(-k),.}$, where $a_{(-k),.}$ is homogeneous of order $-k$ in the $\xi$ variable. In this case $(x,\xi)\mapsto\left|\xi\right|^{k}a_{(-k),.}(x,\xi)$ determines a unique function in $C^{\infty}(\mathbb{S}_{+}^{n}\times\mathbb{S}^{n-1})$ given by \[(z,\omega)\in\mathbb{S}_{+}^{n}\times\mathbb{S}^{n-1}\mapsto a_{(-k),.}(RC^{-1}(z),\omega).\]

\vspace{0,3cm}

(ii) $a\sim\sum_{j=0}^{\infty}a_{.,(-j)}$, where $a_{.,(-j)}$ is homogeneous of order $-j$ in the $x$ variable. In this case $(x,\xi)\mapsto\left|x\right|^{j}a_{.,(-j)}(x,\xi)$ determines a unique function in $C^{\infty}(\mathbb{S}^{n-1}\times\mathbb{S}_{+}^{n})$, similarly as above.

\vspace{0,3cm}

(iii) $a_{(-k),.}\sim\sum_{j=0}^{\infty}a_{(-k),(-j)},\,\,\,\, a_{.,(-j)}\sim\sum_{k=0}^{\infty}a_{(-k),(-j)}$, where $a_{(-k),(-j)}$ is homogeneous of order $-k$ in the $\xi$ variable and homogeneous of order $-j$ in the $x$ variable. In this case $(x,\xi)\mapsto\left|x\right|^{j}\left|\xi\right|^{k}a_{(-k),(-j)}(x,\xi)$ determines a unique function in $C^{\infty}(\mathbb{S}^{n-1}\times\mathbb{S}^{n-1})$, similarly as above.

\vspace{0,3cm}

Conversely, any smooth function that has such an asymptotic expansion defines a classical SG-symbol, as can be shown using Borel Theorem.

In particular, $a_{(0),.}$ can be identified to a function in $C^{\infty}(\mathbb{S}_{+}^{n}\times\mathbb{S}^{n-1})$ and $a_{.,(0)}$ can be identified to a function in $C^{\infty}(\mathbb{S}^{n-1}\times\mathbb{S}_{+}^{n})$. We can even associate the pair $\left(a_{(0),.},a_{.,(0)}\right)$ in an obvious way to a unique function \[\sigma(a)\in C^{\infty}\left(\partial\left(\mathbb{S}_{+}^{n}\times\mathbb{S}_{+}^{n}\right)\right),\]where $\partial\left(\mathbb{S}_{+}^{n}\times\mathbb{S}_{+}^{n}\right)=\mathbb{S}_{+}^{n}\times\mathbb{S}^{n-1}\cup\mathbb{S}^{n-1}\times\mathbb{S}_{+}^{n}$. We are making the identification $\partial\left(\mathbb{S}_{+}^{n}\right)=\mathbb{S}^{n-1}$. For our proposes, it is convenient to call $\sigma(a)$ the principal symbol of $a$. 

It is not hard to prove that if $\sigma(a)$ vanishes, then $a\in S^{-1,-1}_{cl}(\mathbb{R}^n\times\mathbb{R}^n)$.

We can make the same definitions for symbols that assume values in Banach spaces. If $E$ is a complex Banach space, then we define the spaces $S^{0,0}(\mathbb{R}^{n}\times\mathbb{R}^{n},E)$ and $S_{cl}^{0,0}(\mathbb{R}^{n}\times\mathbb{R}^{n},E)$ in the same way as before, but using $C^{\infty}(\mathbb{S}_{+}^{n}\times\mathbb{S}_{+}^{n},E)$, instead of $C^{\infty}(\mathbb{S}_{+}^{n}\times\mathbb{S}_{+}^{n})$.%
\footnote{Sometimes the concept of \textquotedblleft{}operator-valued symbol\textquotedblright{} involves the action of groups of unitary operators on the symbols \cite{KapanadzeSchulze,Schrohebou,Schulze}. In our definition of $S_{cl}^{0,0}(\mathbb{R}^{n}\times\mathbb{R}^{n},E)$, however, we do not consider any action. Hence, for instance, if $F$ and $G$ are Banach spaces, $a\in S^{0,0}(\mathbb{R}^{n}\times\mathbb{R}^{n},\mathcal{B}(F,G))$
if, and only if, $\left\Vert \partial_{x}^{\beta}\partial_{\xi}^{\alpha}a(x,\xi)\right\Vert _{\mathcal{B}(F,G)}\le C\left\langle x\right\rangle^{-|\beta|}\left\langle\xi\right\rangle^{-|\alpha|}$.
}

\selectlanguage{english}%
To each symbol $S^{\mu,\nu}(\mathbb{R}^{n}\times\mathbb{R}^{n})$
and $S_{cl}^{\mu,\nu}(\mathbb{R}^{n}\times\mathbb{R}^{n})$ we can
define a pseudodifferential operator.
\begin{defn}
Let $a\in S^{\mu,\nu}(\mathbb{R}^{n}\times\mathbb{R}^{n})$.
A pseudodifferential operator $A=op(a):\mathcal{S}(\mathbb{R}^{n})\to\mathcal{S}(\mathbb{R}^{n})$ with symbol $a$ is an operator such that, for all $u\in\mathcal{S}(\mathbb{R}^{n})$,
the function $Au$ is the function which assumes at each point $x\in\mathbb{R}^{n}$
the value \[
Au(x)=\frac{1}{\left(2\pi\right)^{n}}\int e^{ix.\xi}a(x,\xi)\hat{u}(\xi)d\xi,\]
where $\hat{u}(\xi)=\int e^{-ix\xi}u(x)dx.$
\end{defn}
The assignment of a symbol to its operator is injective. This allows us to identify
a symbol with its operator. Hence we call classical
those operators whose symbols are classical. 

As in the usual pseudodifferential theory, the zero order operators, that is, the ones whose symbols $a\in S_{cl}^{\mu\nu}(\mathbb{R}^{n}\times\mathbb{R}^{n})$ are such that $\mu\le0$ and $\nu\le0$,
extend continuously to bounded operators $A=op(a):L^{2}(\mathbb{R}^{n})\to L^{2}(\mathbb{R}^{n})$.
Moreover if $\mu<0$ and $\nu<0$, then
$A$ is a compact operator \cite[Theorem 1.4.2]{Rodinonicola}.

The following theorem is an immediate consequence of \cite[Proposition 3.2]{Nicolakteoria}.

\begin{thm}
\label{thm:norma modulo compactos}  Let $A=op(a)$, $a\in S_{cl}^{0,0}(\mathbb{R}^{n}\times\mathbb{R}^{n})$.
Then \[
\begin{array}{c}
\inf_{C\in\mathcal{K}(L^{2}(\mathbb{R}^{n}))}\left\Vert A+C\right\Vert _{\mathcal{B}(L^{2}(\mathbb{R}^{n}))}=\max\left\{ \sup_{x\in\mathbb{R}^{n},|\xi|=1}\left|a_{(0),.}(x,\xi)\right|,\sup_{\xi\in\mathbb{R}^{n},|x|=1}\left|a_{.,(0)}(x,\xi)\right|\right\} ,\end{array}\]

\end{thm}

As a consequence, we conclude that if $a\in S^{0,0}_{cl}(\mathbb {R}^n\times \mathbb {R}^n)$, then the operator $A=op(a):L^{2}(\mathbb{R}^{n})\to L^{2}(\mathbb{R}^{n})$ is Fredholm if, and only if, $\sigma(a)$ never vanishes.

\section{Boutet de Monvel calculus in the classical SG-framework.}

In this section we recall the main properties of the Boutet de Monvel
calculus with classical SG-symbols. A full description of this calculus
can be found in the work of Schrohe \cite{SchroheSGboutet}
for symbols that are not necessarily classical and in the book and article of
Kapanadze and Schulze \cite{KapanadzeSchulze,KapanadzeSchulzearticle}.

We use the following notation: if $x\in\mathbb{R}^{n}$,
then $x=(x',x_{n})$, where $x'\in\mathbb{R}^{n-1}$ and $x_{n}\in\mathbb{R}$.
Similarly for $\xi\in\mathbb{R}^{n}$.

We are interested in the symbols that satisfy the so-called transmission property, in the nomenclature of Rempel and Schulze \cite{RempelSchulze}. We give this definition only for zero-order operators. In what follows $\mathcal{F}$ denotes the Fourier transform $\mathcal{F}(\varphi)=\int e^{ix_{n}\xi_{n}}\varphi(x_{n})dx_{n}$,  $r^{+}:\mathcal{D}'(\mathbb{R}^{n})\to\mathcal{D}'(\mathbb{R}_{+}^{n})$ the restriction operator,  $e^{+}:\mathcal{S}(\mathbb{R}_{+}^{n})\to L^{2}(\mathbb{R}^{n})$ the extension operator to zero on $\mathbb{R}^{n}\backslash\mathbb{R}_{+}^{n}$ and $e^{-}:\mathcal{S}(\mathbb{R}_{-}^{n})\to L^{2}(\mathbb{R}^{n})$ the extension operator to zero on $\overline{\mathbb{R}_{+}^{n}}$,  where $\mathbb{R}_{-}^{n}=\left\{ x\in\mathbb{R}^{n};\, x_{n}<0\right\}$. Finally we denote by $\mathcal{H}_{0}$ \label{H} the set $\mathcal{F}\left(e^{+}\left(\mathcal{S}(\mathbb{R}_{+})\right)\right)\oplus\mathcal{F}\left(e^{-}\left(\mathcal{S}(\mathbb{R}_{-})\right)\right)\oplus\mathbb{C}$.

We say that a symbol of order $(0,0)$, $p\in S_{cl}^{0,0}(\mathbb{R}^{n}\times\mathbb{R}^{n})$,
satisfies the transmission property
if for all $k\in\mathbb{N}_{0}$, we have \[
\left((x',\xi',\xi_{n})\mapsto \left(\partial_{x_{n}}^{k}p\right)(x',0,\xi',\left\langle \xi' \right\rangle \xi_{n})\right)\in S^{0,-k}(\mathbb{R}_{x'}^{n-1}\times\mathbb{R}_{\xi'}^{n-1})\hat{\otimes}\mathcal{H}_{0,\xi_{n}},\]
where $\hat{\otimes}$ denotes the completed projective tensor product.
 
We denote the class of symbols that
satisfy the transmission property by $S_{cl}^{0,0}(\mathbb{R}^{n}\times\mathbb{R}^{n})_{tr}$.

As an example, let $\tilde{a}\in S_{cl}^{0,0}(\mathbb{R}^{n-1}\times\mathbb{R}^{n-1})$,
$\chi$ and $\psi$ be two functions in $C_{c}^{\infty}(\mathbb{R})$.
We define $a(x,\xi):=\chi\left(\frac{x_{n}}{\left\langle x' \right\rangle}\right)\psi\left(\frac{\xi_{n}}{\left\langle \xi' \right\rangle}\right)\tilde{a}(x',\xi')$.
Hence $a\in S_{cl}^{0,0}(\mathbb{R}^{n}\times\mathbb{R}^{n})_{tr}$.

The transmission property has a very important implication. If $a\in S_{cl}^{\mu,\nu}(\mathbb{R}^{n}\times\mathbb{R}^{n})_{tr}$ and $u\in\mathcal{S}(\mathbb{R}_{+}^{n})$, then $r^{+}op(a)e^{+}\left(u\right)\in\mathcal{S}(\mathbb{R}_{+}^{n})$: that is, not only it is smooth on $\mathbb{R}_{+}^{n}$, which is true as the operator $op(a)$ is pseudo-local, but also the function has a smooth extension to a Schwartz function on all $\mathbb{R}^{n}$. Hence $r^{+}op(a)e^{+}$ is an operator from $\mathcal{S}(\mathbb{R}_{+}^{n})$ to $\mathcal{S}(\mathbb{R}_{+}^{n})$.

A classical SG-Boutet de Monvel operator is an operator of the form

\[
\left(\begin{array}{cc}
P_{+}+G & K\\
T & S\end{array}\right):\begin{array}{c}
\mathcal{S}(\mathbb{R}_{+}^{n})\\
\oplus\\
\mathcal{S}(\mathbb{R}^{n-1})\end{array}\to\begin{array}{c}
\mathcal{S}(\mathbb{R}_{+}^{n})\\
\oplus\\
\mathcal{S}(\mathbb{R}^{n-1})\end{array},\]
where $P_{+}=r^{+}Pe^{+}$, with $P$ a classical SG-pseudodifferential
operator that satisfies the transmission property. $T$ is called a trace operator. The class of all trace operators contains the trace operators of classical boundary value problems. Finally the other terms make this into an algebra closed under
adjoints, for zero order operators, and that contains the parametrices of the elliptic operators. The entry $P_{+}$ is called the pseudodifferential
part, $G$ is called a singular Green operator, $K$ is called a Poisson
operator, $T$ is called a Trace operator and $S$ is called the pseudodifferential
operator of the border.
\footnote{Actually, we could define operators from $\mathcal{S}(\mathbb{R}_{+}^{n})^{N}\oplus\mathcal{S}(\mathbb{R}^{n-1})^{M}$ to $\begin{array}{c}
\mathcal{S}(\mathbb{R}_{+}^{n})^{N'}\oplus\mathcal{S}(\mathbb{R}^{n-1})^{M'}\end{array}$, for all $N$, $N$', $M$ and $M$ that belong to $\mathbb{N}_{0}$.
Nevertheless, for the study of the $K$-Theory only the case $N=N'=M=M'=1$
matters. (see \cite[Section 1.5]{MeloSchrohe})%
}

Using order reducing operators, see for instance \cite{SchroheSGboutet}, for index proposes we can restrict the study to Boutet de Monvel Operators of order $(0,0)$ and type zero. We give precise definitions only in the case where the type is zero.

Our notation is similar to Grubb's
\cite{Grubbamarelo,Grubbverde}. We will always denote by $\mathcal{S}_{+}$
and $\mathcal{S}_{++}$ the spaces $\mathcal{S}(\mathbb{R}_{+})$
and $\mathcal{S}(\mathbb{R}_{+}\times\mathbb{R}_{+})$, respectively. Without changing notation, we now consider RC also defined on $\mathbb{R}^{n-1}$ with values in $\mathbb{S}^{n-1}_{+}$. Thus we regard $\mathbb{S}^{n-1}_{+}\times\mathbb{S}^{n-1}_{+}$ as a compactification of $\mathbb{R}^{n-1}\times\mathbb{R}^{n-1}$.

\begin{defn}
The space $S_{cl}^{0,0}(\mathbb{R}^{n-1},\mathbb{R}^{n-1},\mathcal{S}_{+})$
is the set of functions $f\in C^{\infty}(\mathbb{R}^{n-1}\times\mathbb{R}_{+}\times\mathbb{R}^{n-1})$
such that

1) For each fixed $(x',\xi')$ the function $\left(x_{n}\mapsto f(x',x_{n},\xi')\right)\in\mathcal{S}_{+}$.

2) Let $f^{[0]}\in C^{\infty}(\mathbb{R}^{n-1}\times\mathbb{R}_{+}\times\mathbb{R}^{n-1})$
be the function $f^{[0]}(x',x_{n},\xi')=\left\langle \xi' \right\rangle^{-\frac{1}{2}}f(x',\left\langle\xi'\right\rangle^{-1}x_{n},\xi')$.
Then the function $f^{[0]}\circ\left(RC^{-1}\times Id\times RC^{-1}\right):\left(\mathbb{S}_{+}^{n-1}\cap\mathbb{R}_{+}^{n}\right)\times\left(\mathbb{S}_{+}^{n-1}\cap\mathbb{R}_{+}^{n}\right)\to\mathcal{S}_{+}$
given by \[
\left(z,w\right)\mapsto\left(x_{n}\mapsto f\left(RC^{-1}(z),x_{n},RC^{-1}(w)\right)\right)\]
 can be extended to a function in $C^{\infty}(\mathbb{S}^{n-1}\times\mathbb{S}^{n-1},\mathcal{S}_{+})$. 

Let $(\mu,\nu)\in\mathbb{R}^{2}$. The space $S_{cl}^{\mu,\nu}(\mathbb{R}^{n-1},\mathbb{R}^{n-1},\mathcal{S}_{+})$
consists of functions $f\in C^{\infty}(\mathbb{R}^{n-1}\times\mathbb{R}{}_{+}\times\mathbb{R}^{n-1})$
such that \[
\left((x',x_{n},\xi')\mapsto\left\langle x' \right\rangle^{-\nu}\left\langle \xi' \right\rangle^{-\mu}f(x',x_{n},\xi')\right)\in S_{cl}^{0,0}(\mathbb{R}^{n-1},\mathbb{R}^{n-1},\mathcal{S}_{+}).\]

\end{defn}

\begin{defn}
The space $S_{cl}^{0,0}(\mathbb{R}^{n-1},\mathbb{R}^{n-1},\mathcal{S}_{++})$
consists of functions $g\in C^{\infty}\left(\mathbb{R}^{n-1}\times\mathbb{R}{}_{++}\times\mathbb{R}^{n-1}\right)$
such that

1) For each fixed $(x',\xi')$ the function $\left(\left(x_{n},y_{n}\right)\mapsto g(x',x_{n},y_{n},\xi')\right)\in\mathcal{S}_{++}$.

2) Let $g^{[0]}\in C^{\infty}(\mathbb{R}^{n-1}\times\mathbb{R}{}_{++}\times\mathbb{R}^{n-1})$
be the function \[g^{[0]}(x',x_{n},y_{n},\xi')=\left\langle \xi' \right\rangle^{-1}g(x',\left\langle \xi' \right\rangle^{-1}x_{n},\left\langle \xi' \right\rangle^{-1}y_{n},\xi').\]
Then the function $g^{[0]}\circ\left(RC^{-1}\times Id\times RC^{-1}\right):\left(\mathbb{S}_{+}^{n-1}\cap\mathbb{R}_{+}^{n}\right)\times\left(\mathbb{S}_{+}^{n-1}\cap\mathbb{R}_{+}^{n}\right)\to\mathcal{S}_{++}$
given by \[
\left(z,w\right)\mapsto\left(\left(x_{n},y_{n}\right)\mapsto g\left(RC^{-1}(z),x_{n},y_{n},RC^{-1}(w)\right)\right)\]
 can be extended to a function in $C^{\infty}(\mathbb{S}^{n-1}\times\mathbb{S}^{n-1},\mathcal{S}_{++})$. 

Let $(\mu,\nu)\in\mathbb{R}^{2}$. The space $S_{cl}^{\mu,\nu}(\mathbb{R}^{n-1},\mathbb{R}^{n-1},\mathcal{S}_{++})$
consists of functions $g\in C^{\infty}(\mathbb{R}^{n-1}\times\mathbb{R}{}_{++}\times\mathbb{R}^{n-1})$ such that
 \[
\left((x',x_{n},y_{n},\xi')\mapsto\left\langle x' \right\rangle^{-\nu}\left\langle \xi' \right\rangle^{-\mu}g(x',x_{n},y_{n},\xi')\right)\in S_{cl}^{0,0}(\mathbb{R}^{n-1},\mathbb{R}^{n-1},\mathcal{S}_{++}).\]

\end{defn}

\begin{defn}
$(i)$ A classical SG-trace operator $op(t):\mathcal{S}(\mathbb{R}_{+}^{n})\to\mathcal{S}(\mathbb{R}^{n-1})$
of order $(\mu,\nu)\in\mathbb{Z}^{2}$ and type $0$
is an operator of the form \[
op(t)u(x')=\frac{1}{(2\pi)^{n-1}}\int_{\mathbb{R}^{n-1}}e^{ix'\xi'}\int_{0}^{\infty}t(x',x_{n},\xi')\acute{u}(\xi',x_{n})dx_{n}d\xi',\]
where $t\in S_{cl}^{\mu,\nu}(\mathbb{R}^{n-1},\mathbb{R}^{n-1},\mathcal{S}_{+})$
and $\acute{u}(\xi',x_{n})=\int e^{-ix'\xi'}u(x',x_{n})dx'$. The
set of these operators is denoted by $\mathcal{T}^{(\mu,\nu),0}(\mathbb{R}_{+}^{n})$.

$(ii)$ A classical SG-Poisson operator $op(k):\mathcal{S}(\mathbb{R}^{n-1})\to\mathcal{S}(\mathbb{R}_{+}^{n})$
of order $(\mu,\nu)\in\mathbb{Z}^{2}$ is an operator defined by \[
op(k)u(x',x_{n})=\frac{1}{(2\pi)^{n-1}}\int_{\mathbb{R}^{n-1}}e^{ix'\xi'}k(x',x_{n},\xi')\hat{u}(\xi')d\xi',\]
where $k\in S_{cl}^{\mu,\nu}(\mathbb{R}^{n-1},\mathbb{R}^{n-1},\mathcal{S}_{+})$
and $\hat{u}(\xi')=\int e^{-ix'\xi'}u(x')dx'$. The set of these operators
is denoted by $\mathcal{K}^{(\mu,\nu)}(\mathbb{R}_{+}^{n})$.

$(iii)$ A classical SG-singular Green operator $op(g):\mathcal{S}(\mathbb{R}^{n}_{+})\to\mathcal{S}(\mathbb{R}_{+}^{n})$
of order $(\mu,\nu)\in\mathbb{Z}^{2}$ and type $0$
is an operator of the form\[
op(g)u(x)=\frac{1}{(2\pi)^{n-1}}\int_{\mathbb{R}^{n-1}}e^{ix'\xi'}\int_{0}^{\infty}g(x',x_{n},y_{n},\xi')\acute{u}(\xi',y_{n})dy_{n}d\xi',\]
where $g\in S_{cl}^{\mu,\nu}(\mathbb{R}^{n-1},\mathbb{R}^{n-1},\mathcal{S}_{++})$
and $\acute{u}(\xi',x_{n})=\int e^{-ix'\xi'}u(x',x_{n})dx'$. The
set of these operators is denoted by $\mathcal{G}^{(\mu,\nu),0}(\mathbb{R}_{+}^{n})$.

$(iv)$ A classical SG-pseudodifferential operator on the boundary $op(s):\mathcal{S}(\mathbb{R}^{n-1})\to\mathcal{S}(\mathbb{R}^{n-1})$
of order $(\mu,\nu)\in\mathbb{Z}^{2}$ is a pseudodifferential operator with symbol $s\in S_{cl}^{\mu,\nu}(\mathbb{R}^{n-1}\times\mathbb{R}^{n-1})$.
This means that\[
op(s)u(x')=\frac{1}{(2\pi)^{n-1}}\int_{\mathbb{R}^{n-1}} e^{ix'\xi'}s(x',\xi')\hat{u}(\xi')d\xi',\]
where as usual $\hat{u}(\xi'):=\int e^{-ix'\xi'}u(x')dx'$.

$(v)$ The pseudodifferential part $P_{+}:\mathcal{S}(\mathbb{R}_{+}^{n})\to\mathcal{S}(\mathbb{R}_{+}^{n})$ of a Boutet de Monvel operator of order $(0,0)$
acting on $\mathbb{R}_{+}^{n}$
is an operator of the form $P_{+}=r^{+}op(p)e^{+}$, where $p$
is a symbol that belongs to $S_{cl}^{0,0}(\mathbb{R}^{n}\times\mathbb{R}^{n})_{tr}$,
that is, it satisfies the transmission property.
\end{defn}

\begin{thm}
\cite[Section 4]{SchroheSGboutet} The operators above have the following continuous extensions:

(i) If $p\in S_{cl}^{0,0}(\mathbb{R}^{n}\times\mathbb{R}^{n})_{tr}$,
then $P_{+}$ has a continuous extension $P_{+}:L^{2}(\mathbb{R}_{+}^{n})\to L^{2}(\mathbb{R}_{+}^{n})$.

(ii) If $G\in\mathcal{G}^{(0,0),0}(\mathbb{R}_{+}^{n})$, then $G$ has a continuous extension $G:L^{2}(\mathbb{R}_{+}^{n})\to L^{2}(\mathbb{R}_{+}^{n})$.

(iii) If $T\in\mathcal{T}^{(0,0),0}(\mathbb{R}_{+}^{n})$, then $T$ has a continuous extension $T:L^{2}(\mathbb{R}_{+}^{n})\to L^{2}(\mathbb{R}^{n-1})$.

(iv) If $K\in\mathcal{K}^{(0,0)}(\mathbb{R}_{+}^{n})$, then $K$ has a continuous extension $K:L^{2}(\mathbb{R}^{n-1})\to L^{2}(\mathbb{R}_{+}^{n})$.
\end{thm}
In order to understand the regularizing operators, we use the
following definitions and propositions - see \cite{SchroheSGboutet}.
\begin{defn}
1) We say that the pseudodifferential part $P_{+}$ is regularizing
if $P=op(p)$, where \[
p\in\cap_{\mu,\nu\in\mathbb{Z}}S_{cl}^{\mu,\nu}(\mathbb{R}^{n}\times\mathbb{R}^{n})=S^{-\infty,-\infty}(\mathbb{R}^{n}\times\mathbb{R}^{n}).\]

$P_{+}:\mathcal{S}(\mathbb{R}_{+}^{n})\to\mathcal{S}(\mathbb{R}_{+}^{n})$
is regularizing iff it is an integral operator whose kernel belongs to $\mathcal{S}(\mathbb{R}_{+}^{n}\times\mathbb{R}_{+}^{n})$.

2) We say that a singular Green operator $G$ is regularizing of type
$0$ if $G=op(g)$, where \[
g\in\cap_{\mu,\nu\in\mathbb{Z}}S_{cl}^{\mu,\nu}(\mathbb{R}^{n-1},\mathbb{R}^{n-1},\mathcal{S}_{++}).\]

A singular Green operator of type zero $G:\mathcal{S}(\mathbb{R}_{+}^{n})\to\mathcal{S}(\mathbb{R}_{+}^{n})$
is regularizing iff it is an integral operator whose kernel belongs to $\mathcal{S}(\mathbb{R}_{+}^{n}\times\mathbb{R}_{+}^{n})$.
The regularizing singular Green operators of type zero are therefore
equal to the regularizing pseudodifferential operators $P_{+}$.

3) We say that a trace operator $T$ is regularizing of type $0$
if $T=op(t)$, where \[
t\in\cap_{\mu,\nu\in\mathbb{Z}}S_{cl}^{\mu,\nu}(\mathbb{R}^{n-1},\mathbb{R}^{n-1},\mathcal{S}_{+}).\]

A trace operator of type zero $T:\mathcal{S}(\mathbb{R}_{+}^{n})\to\mathcal{S}(\mathbb{R}^{n-1})$
is regularizing iff it is an integral operator whose kernel belongs to $\mathcal{S}(\mathbb{R}^{n-1}\times\mathbb{R}_{+}^{n})$.

4) We say that a Poisson operator $K$ is regularizing if $K=op(k)$,
where \[
k\in\cap_{\mu,\nu\in\mathbb{Z}}S_{cl}^{\mu,\nu}(\mathbb{R}^{n-1},\mathbb{R}^{n-1},\mathcal{S}_{+}).\]

A Poisson operator $K:\mathcal{S}(\mathbb{R}^{n-1})\to\mathcal{S}(\mathbb{R}_{+}^{n})$
is regularizing iff it is an integral operator whose kernel belongs to $\mathcal{S}(\mathbb{R}_{+}^{n}\times\mathbb{R}^{n-1})$.
\end{defn}
The following theorem says that away from the boundary $\{x\in\overline{\mathbb{R}^{n}_{+}};x_{n}=0\}$, the singular Green, the trace and the Poisson operators are regularizing. 
\begin{thm}
\label{thm:operadores regularizantes de boutet}
\cite[Theorem 2.10]{SchroheSGboutet} Let $\phi\in C_{c}^{\infty}(\mathbb{R})$
be a function that assumes the value $1$ in a neighborhood of zero.
Let us define the function $\Phi\in C^{\infty}(\mathbb{R}^{n})$
by \[
\Phi(x):=1-\phi\left(\frac{x_{n}}{\left\langle x' \right\rangle}\right).\]

Let $G$, $T$ and $K$ be a singular Green operator, a trace operator
and a Poisson operator of order $(0,0)$ and type
$0$ - when a type is defined, respectively. Then $\Phi K$ is a regularizing
Poisson operator, $T\Phi$ is a regularizing trace operator of type
zero, $G\Phi$ and $\Phi G$ are regularizing singular Green operator
of type zero.
\end{thm}

\begin{defn}
\label{def:sgboutet}
The set of classical SG-Boutet de Monvel operators of order $(0,0)$ 
and type $0$, denoted by $\mathcal{B}^{(0,0),0}(\mathbb{R}^{n}_{+})$,
is the set of operators of the form\[
A=\left(\begin{array}{cc}
P_{+}+G & K\\
T & S\end{array}\right):\begin{array}{c}
\mathcal{S}(\mathbb{R}_{+}^{n})\\
\oplus\\
\mathcal{S}(\mathbb{R}^{n-1})\end{array}\to\begin{array}{c}
\mathcal{S}(\mathbb{R}_{+}^{n})\\
\oplus\\
\mathcal{S}(\mathbb{R}^{n-1})\end{array},\]
where $P=r^{+}Op(p)e^{+}$, with $p\in S_{cl}^{0,0}(\mathbb{R}^{n}\times\mathbb{R}^{n})_{tr}$,
$G\in\mathcal{G}^{(0,0),0}(\mathbb{R}_{+}^{n})$, $K\in\mathcal{K}^{(0,0)}(\mathbb{R}_{+}^{n})$,
$T\in\mathcal{T}^{(0,0),0}(\mathbb{R}_{+}^{n})$ and $S=Op(s)$, with
$s\in S_{cl}^{0,0}(\mathbb{R}^{n-1}\times\mathbb{R}^{n-1})$.

We say that $A$ is regularizing of type $0$ if all the terms in
the matrix are regularizing.
\end{defn}
We note that $\mathcal{B}^{(0,0),0}(\mathbb{R}^{n}_{+})$ is closed under
composition and taking formal adjoints. Hence it is an $*$-algebra. Furthermore its elements extend to bounded operators on $L^{2}(\mathbb{R}_{+}^{n})\oplus L^{2}(\mathbb{R}^{n-1})$.

\section{Principal Boundary Value Symbol.}

As in the case of classical SG-symbols, explained in section \ref{sec:sec1}, we say that a function $f$ belongs to $S_{cl}^{0,0}(\mathbb{R}^{n-1},\mathbb{R}^{n-1},\mathcal{S}_{+})$ iff $f^{[0]}$ has a smooth extension to the boundary of the compactification $\mathbb{S}^{n-1}_{+}\times\mathbb{S}^{n-1}_{+}$. We can then take Taylor expansions in the variables $t=\frac{1}{\left|x'\right|}$ and $s=\frac{1}{\left|\xi'\right|}$ (the other variables are $\Omega_{1}=\frac{x'}{\left|x'\right|}$ and $\Omega_{2}=\frac{\xi'}{\left|\xi'\right|}$) to obtain the following asymptotic expansions:

\vspace{0,3cm}

(i) $f^{[0]}\sim\sum_{k=0}^{\infty}f^{[0]}_{(-k),.}$, where $f^{[0]}_{(-k),.}$ is homogeneous of order $-k$ in the $\xi'$ variable. Using the radial compactification, the function \[(x',\xi')\mapsto\left(x_{n}\mapsto\left|\xi'\right|^{k}f^{[0]}_{(-k),.}(x',x_{n},\xi')\right)\] can be identified with the function
\[(z,\omega)\in\mathbb{S}_{+}^{n-1}\times\mathbb{S}^{n-2}\mapsto \left(x_{n}\mapsto f^{[0]}_{(-k),.}(RC^{-1}(z),x_{n},\omega) \right) \] in $C^{\infty}\left(\mathbb{S}_{+}^{n-1}\times\mathbb{S}^{n-2},\mathcal{S}_{+}\right)$.

\vspace{0,3cm}

(ii) $f^{[0]}\sim\sum_{j=0}^{\infty}f^{[0]}_{.,(-j)}$, where $f^{[0]}_{.,(-j)}$ is homogeneous of order $-j$ in the $x'$ variable. Using the radial compactification the function \[(x',\xi')\mapsto\left(x_{n}\mapsto\left|x'\right|^{j}f^{[0]}_{.,(-j)}(x',x_{n},\xi')\right)\] can be identified to a function in $C^{\infty}\left(\mathbb{S}^{n-2}\times\mathbb{S}_{+}^{n-1},\mathcal{S}_{+}\right)$ similarly as in (i).

\vspace{0,3cm}

(iii) $f^{[0]}_{(-k),.}\sim\sum_{j=0}^{\infty}f^{[0]}_{(-k),(-j)},\,\,\,\, f^{[0]}_{.,(-j)}\sim\sum_{k=0}^{\infty}f^{[0]}_{(-k),(-j)}$, where $f^{[0]}_{(-k),(-j)}$ is homogeneous of order $-k$ in the $\xi'$ variable and homogeneous of order $-j$ in the $x'$ variable. Using the radial compactification the function \[(x',\xi')\mapsto\left(x_{n}\mapsto\left|x'\right|^{j}\left|\xi'\right|^{k}f^{[0]}_{(-k),(-j)}(x',x_{n},\xi')\right)\] can be identified to a function in $C^{\infty}\left(\mathbb{S}^{n-2}\times\mathbb{S}^{n-2},\mathcal{S}_{+}\right)$ similarly as in (i).

\vspace{0,3cm}

In particular, $f^{[0]}_{(0),.}$ can be identified to a function in $C^{\infty}(\mathbb{S}_{+}^{n-1}\times\mathbb{S}^{n-2},\mathcal{S}_{+})$ and $f^{[0]}_{.,(0)}$ can be identified to a function in $C^{\infty}(\mathbb{S}^{n-2}\times\mathbb{S}_{+}^{n-1},\mathcal{S}_{+})$. We can even associate the pair $\left(f^{[0]}_{(0),.},f^{[0]}_{.,(0)}\right)$ in a obvious way to a unique function in $C^{\infty}\left(\partial\left(\mathbb{S}_{+}^{n-1}\times\mathbb{S}_{+}^{n-1}\right),\mathcal{S}_{+}\right)$, where $\partial\left(\mathbb{S}_{+}^{n-1}\times\mathbb{S}_{+}^{n-1}\right)=\mathbb{S}_{+}^{n-1}\times\mathbb{S}^{n-2}\cup\mathbb{S}^{n-2}\times\mathbb{S}_{+}^{n-1}$. We are making the identification $\partial\left(\mathbb{S}_{+}^{n-1}\right)=\mathbb{S}^{n-2}$. We consider two complex-valued functions $f_{(0),.}$ and $f_{.,(0)}$ defined on $\mathbb{R}^{n-1}\times\mathbb{R}_{+}\times\left(\mathbb{R}^{n-1}\backslash\{0\}\right)$ and $\left(\mathbb{R}^{n-1}\backslash\{0\}\right)\times\mathbb{R}_{+}\times\mathbb{R}^{n-1}$, respectively, by 
\begin{equation}f_{(0),.}\left(x',x_{n},\xi'\right):=\lim_{\lambda\to\infty}\lambda^{-\frac{1}{2}}f\left(x',\lambda^{-1}x_{n},\lambda\xi'\right)=\left|\xi'\right|^{\frac{1}{2}}f_{(0),.}^{[0]}\left(x',\left|\xi'\right|x_{n},\xi'\right)\label{eq:f0,.}\end{equation} and \label{f0,.}\begin{equation}f_{.,(0)}\left(x',x_{n},\xi'\right):=\lim_{\lambda\to\infty}f\left(\lambda x',x_{n},\xi'\right)=\left\langle \xi'\right\rangle ^{\frac{1}{2}}f_{.,(0)}^{[0]}\left(x',\left\langle \xi'\right\rangle x_{n},\xi'\right)\label{eq:f.,0}.\end{equation}

In precisely the same way, we define $g_{(0),.}$ and $g_{.,(0)}$ for a function $g\in S_{cl}^{0,0}(\mathbb{R}^{n-1},\mathbb{R}^{n-1},\mathcal{S}_{++})$. We only need to replace $S_{+}$ by $S_{++}$, $\left\langle \xi'\right \rangle^{\frac{1}{2}}$ by $\left\langle \xi' \right\rangle$ and  $\left| \xi'\right |^{\frac{1}{2}}$ by $\left| \xi' \right|$.

\begin{defn}
Let $p\in S_{cl}^{0,0}(\mathbb{R}^{n}\times\mathbb{R}^{n})_{tr}$, $g\in S_{cl}^{0,0}(\mathbb{R}^{n-1},\mathbb{R}^{n-1},\mathcal{S}_{++})$, $k\in S_{cl}^{0,0}(\mathbb{R}^{n-1},\mathbb{R}^{n-1},\mathcal{S}_{+})$, $t\in S_{cl}^{0,0}(\mathbb{R}^{n-1},\mathbb{R}^{n-1},\mathcal{S}_{+})$ and $s\in S_{cl}^{0,0}(\mathbb{R}^{n-1}\times\mathbb{R}^{n-1})$. Then to each classical SG-Boutet de Monvel operator \[A:=\left(\begin{array}{cc}
op(p)_{+}+op(g) & op(k)\\
op(t) & op(s)\end{array}\right):\begin{array}{c}
\mathcal{S}(\mathbb{R}_{+}^{n})\\
\oplus\\
\mathcal{S}(\mathbb{R}^{n-1})\end{array}\to\begin{array}{c}
\mathcal{S}(\mathbb{R}_{+}^{n})\\
\oplus\\
\mathcal{S}(\mathbb{R}^{n-1})\end{array},\]we assign two functions: $A_{(0),.}:\mathbb{R}^{n-1}\times\left(\mathbb{R}^{n-1}\backslash\{0\}\right)\to\mathcal{B}\left(L^{2}(\mathbb{R}_{+})\oplus\mathbb{C}\right)$ and $A_{.,(0)}:\left(\mathbb{R}^{n-1}\backslash\{0\}\right)\times\mathbb{R}^{n-1}\to\mathcal{B}\left(L^{2}(\mathbb{R}_{+})\oplus\mathbb{C}\right)$ given by \[A_{(0),.}(x',\xi'):=\left(\begin{array}{cc}
p_{(0),.}(x',\xi',D)_{+}+g_{(0),.}(x',\xi',D) & k_{(0),.}(x',\xi',D)\\
t_{(0),.}(x',\xi',D) & s_{(0),.}(x',\xi',D)\end{array}\right)\]and\[A_{.,(0)}(x',\xi'):=\left(\begin{array}{cc}
p_{.,(0)}(x',\xi',D)_{+}+g_{.,(0)}(x',\xi',D) & k_{.,(0)}(x',\xi',D)\\
t_{.,(0)}(x',\xi',D) & s_{.,(0)}(x',\xi',D)\end{array}\right),\]where:

(i) the operator $p_{(0),.}(x',\xi',D)_{+}:L^{2}(\mathbb{R}_{+})\to L^{2}(\mathbb{R}_{+})$ is defined by \[p_{(0),.}(x',\xi',D)_{+}u(x_{n})=\frac{1}{2\pi}\int_{-\infty}^{\infty} e^{ix_{n}\xi_{n}}p_{(0),.}(x',\xi',\xi_{n})\widehat{e^{+}u}(\xi_{n})d\xi_{n},\] where $\widehat{e^{+}u}(\xi_{n})=\int_{0}^{\infty}e^{-ix_{n}\xi_{n}}u(x_{n})dx_{n}$;

(ii) the operator $g_{(0),.}(x',\xi',D):L^{2}(\mathbb{R}_{+})\to L^{2}(\mathbb{R}_{+})$ is defined by \[g_{(0),.}(x',\xi',D)u(x_{n})=\int_{0}^{\infty}g_{(0),.}(x',x_{n},y_{n},\xi')u(y_{n})dy_{n};\]

(iii) the operator $k_{(0),.}(x',\xi',D):\mathbb{C}\to L^{2}(\mathbb{R}_{+})$ is defined by \[k_{(0),.}(x',\xi',D)c=k_{(0),.}(x',x_{n},\xi')c;\]

(iv) the operator $t_{(0),.}(x',\xi',D):L^{2}(\mathbb{R}_{+})\to\mathbb{C}$ is defined by \[t_{(0),.}(x',\xi',D)u(x_{n})=\int_{0}^{\infty} t_{(0),.}(x',x_{n},\xi')u(x_{n})dx_{n};\]

(v) the operator $s_{(0),.}(x',\xi',D):\mathbb{C}\to\mathbb{C}$ is defined by \[s_{(0),.}(x',\xi',D)c=s_{(0),.}(x',\xi')c.\]

For $p_{.,(0)}(x',\xi',D)_{+}$, $g_{.,(0)}(x',\xi',D)$, $t_{.,(0)}(x',\xi',D)$, $k_{.,(0)}(x',\xi',D)$ and $s_{.,(0)}(x',\xi',D)$ hold similar definitions.
\end{defn}
The assignment $A\mapsto \left( A_{(0),.},A_{.,(0)} \right)$ is a $*$-homomorphism.

Using the arguments of Rempel and Schulze \cite[Section 2.3.4]{RempelSchulze},
it is straightforward to prove the following theorem, analogous to our Theorem 3. It is important to remark that if $P_{+}+G=P'_{+}+G'$ then the principal symbols of $P$ and $P'$ are equal for $x_{n}>0$ (details can be found in \cite{Tese}).

\begin{thm}
Let $A\in\mathcal{B}^{(0,0),0}(\mathbb{R}^{n}_{+})$ be given by \[
A=\left(\begin{array}{cc}
P_{+}+G & K\\
T & S\end{array}\right),\]
Let $A_{(0),.}$
and $A_{.,(0)}$ be functions defined as above. Let $p_{(0),.}$
and $p_{.,(0)}$ be the first terms of the asymptotic expansion of the symbol of $P$. Then the following estimate holds:\foreignlanguage{brazil}{\[
\inf_{C\in\mathcal{K}(L^{2}(\mathbb{R}_{+}^{n})\oplus L^{2}(\mathbb{R}^{n-1}))}\left\Vert A+C\right\Vert _{\mathcal{B}(L^{2}(\mathbb{R}_{+}^{n})\oplus L^{2}(\mathbb{R}^{n-1}))}=\max\left\{ \sup_{x\in\mathbb{R}_{+}^{n},|\xi|=1}\left|p_{(0),.}(x,\xi)\right|,\sup_{\xi\in\mathbb{R}^{n},x\in\mathbb{R}_{+}^{n}\, ,|x|=1}\left|p_{.,(0)}(x,\xi)\right|,\right.\]
\[
\left.\sup_{x'\in\mathbb{R}^{n-1},|\xi'|=1}\left\Vert A_{(0),.}(x',\xi')\right\Vert _{\mathcal{B}(L^{2}(\mathbb{R}_{+})\oplus\mathbb{C})},\sup_{\xi'\in\mathbb{R}^{n-1},|x'|=1}\left\Vert A_{.,(0)}(x',\xi')\right\Vert _{\mathcal{B}(L^{2}(\mathbb{R}_{+})\oplus\mathbb{C})}\right\} .\]
}
\end{thm}

Let us show now how to assign to each pair $A_{(0),.}$ and $A_{.,(0)}$ a unique function that belongs to $C\left(\partial\left(\mathbb{S}_{+}^{n-1}\times\mathbb{S}_{+}^{n-1}\right),\mathcal{B}\left(L^{2}(\mathbb{R}_{+})\oplus\mathbb{C}\right)\right)$, where $\partial\left(\mathbb{S}_{+}^{n-1}\times\mathbb{S}_{+}^{n-1}\right)=\mathbb{S}_{+}^{n-1}\times\mathbb{S}^{n-2}\cup\mathbb{S}^{n-2}\times\mathbb{S}_{+}^{n-1}$.

We start with the Poisson operator. The others are treated very similarly. Let us consider a function $k\in\mathcal{S}_{cl}^{0,0}\left(\mathbb{R}^{n-1},\mathbb{R}^{n-1},\mathcal{S}_{+}\right)$ and let $op(k):\mathcal{S}(\mathbb{R}^{n-1})\to\mathcal{S}(\mathbb{R}_{+}^{n})$ be a Poisson operator. 

The functions $k_{(0),.}$ and $k_{.,(0)}$, see equations (\ref{eq:f0,.}) and (\ref{eq:f.,0}), define, for each $\left(x',\xi'\right)$, functions $k_{(0),.}(x',\xi',D):\mathbb{C}\to L^{2}(\mathbb{R}_{+})$ and $k_{.,(0)}(x',\xi',D):\mathbb{C}\to L^{2}(\mathbb{R}_{+})$, as we have already seen. Let us now define the unitary operator $\kappa_{\lambda}:L^{2}(\mathbb{R}_{+})\to L^{2}(\mathbb{R}_{+})$ by $\kappa_{\lambda}u(x)=\lambda^{\frac{1}{2}}u(\lambda x)$. Hence it is easy to check that\[\kappa_{\left|\xi'\right|^{-1}}k_{(0),.}(x',\xi',D)c=\left|\xi'\right|^{-\frac{1}{2}}k_{(0),.}(x',\left|\xi'\right|^{-1}x_{n},\xi')c=k_{(0),.}^{[0]}(x',x_{n},\xi')c\]and\[\kappa_{\left\langle \xi' \right\rangle^{-1}}k_{.,(0)}(x',\xi',D)c=\left\langle \xi' \right\rangle^{-\frac{1}{2}}k_{(0),.}(x',\left\langle \xi' \right\rangle^{-1}x_{n},\xi')c=k_{.,(0)}^{[0]}(x',x_{n},\xi')c.\] We know that $k_{(0),.}^{[0]}$ and $k_{.,(0)}^{[0]}$ can be seen as a unique function in $C^{\infty}\left(\partial\left(\mathbb{S}_{+}^{n-1}\times\mathbb{S}_{+}^{n-1}\right),\mathcal{S}_{+}\right)$. As $\mathcal{S}_{+}\hookrightarrow L^{2}(\mathbb{R}_{+})$ continuously and, therefore, $C^{\infty}\left(\partial\left(\mathbb{S}_{+}^{n-1}\times\mathbb{S}_{+}^{n-1}\right),\mathcal{S}_{+}\right)$ has a continuous inclusion in $C^{\infty}\left(\partial\left(\mathbb{S}_{+}^{n-1}\times\mathbb{S}_{+}^{n-1}\right),\mathcal{B}(\mathbb{C},L^{2}(\mathbb{R}_{+}))\right)$, we conclude that $\kappa_{\left|\xi'\right|^{-1}}k_{(0),.}(x',\xi',D)$ and $\kappa_{\left\langle \xi' \right\rangle^{-1}}k_{.,(0)}(x',\xi',D)$ determine in a unique way a function in $C^{\infty}\left(\partial\left(\mathbb{S}_{+}^{n-1}\times\mathbb{S}_{+}^{n-1}\right),\mathcal{B}\left(\mathbb{C},L^{2}(\mathbb{R}_{+})\right)\right)$.

Putting everything together we obtain the following result: \label{unique}

\begin{prop}
\label{pro:protreze}
(i) Let $k\in S_{cl}^{0,0}(\mathbb{R}^{n-1},\mathbb{R}^{n-1},\mathcal{S}_{+})$ and $op(k)$ be a Poisson operator. Let $k_{(0),.}$ and $k_{.,(0)}$ be as defined in equations (\ref{eq:f0,.}) and (\ref{eq:f.,0}). Hence $\kappa_{\left|\xi'\right|^{-1}}k_{(0),.}(x',\xi',D)$ and $\kappa_{\left\langle \xi' \right\rangle^{-1}}k_{.,(0)}(x',\xi',D)$ determine, in a unique way, a function in $C^{\infty}\left(\partial\left(\mathbb{S}_{+}^{n-1}\times\mathbb{S}_{+}^{n-1}\right),\mathcal{B}\left(\mathbb{C},L^{2}(\mathbb{R}_{+})\right)\right)$. Furthermore, the closure of all functions constructed in this way is equal to $C\left(\partial\left(\mathbb{S}_{+}^{n-1}\times\mathbb{S}_{+}^{n-1}\right),\mathcal{B}\left(\mathbb{C},L^{2}(\mathbb{R}_{+})\right)\right)$

\vspace{0,2cm}

(ii) Let $t\in S_{cl}^{0,0}(\mathbb{R}^{n-1},\mathbb{R}^{n-1},\mathcal{S}_{+})$ and $op(t)$ be a trace operator. Let $t_{(0),.}$ and $t_{.,(0)}$ be as defined in equations (\ref{eq:f0,.}) and (\ref{eq:f.,0}). Hence $t_{(0),.}(x',\xi',D)\kappa_{\left|\xi'\right|}$ and $t_{.,(0)}(x',\xi',D)\kappa_{\left\langle \xi' \right\rangle}$ determine, in a unique way, a function in $C^{\infty}\left(\partial\left(\mathbb{S}_{+}^{n-1}\times\mathbb{S}_{+}^{n-1}\right),\mathcal{B}(L^{2}(\mathbb{R}_{+}),\mathbb{C})\right)$. Furthermore the closure of all functions constructed in this way is equal to $C\left(\partial\left(\mathbb{S}_{+}^{n-1}\times\mathbb{S}_{+}^{n-1}\right),\mathcal{B}(L^{2}(\mathbb{R}_{+}),\mathbb{C})\right)$

\vspace{0,2cm}

(iii) Let $g\in S_{cl}^{0,0}(\mathbb{R}^{n-1},\mathbb{R}^{n-1},\mathcal{S}_{++})$ and $op(g)$ be a Green operator. Let $g_{(0),.}$ and $g_{.,(0)}$ be as defined in equations (\ref{eq:f0,.}) and (\ref{eq:f.,0}). Hence $\kappa_{\left|\xi'\right|^{-1}}g_{(0),.}(x',\xi',D)\kappa_{\left|\xi'\right|}$ and $\kappa_{\left\langle \xi' \right\rangle^{-1}}g_{.,(0)}(x',\xi',D)\kappa_{\left\langle \xi' \right\rangle}$ determine, in a unique way, a function in $C^{\infty}\left(\partial\left(\mathbb{S}_{+}^{n-1}\times\mathbb{S}_{+}^{n-1}\right),\mathcal{B}\left(L^{2}(\mathbb{R}_{+})\right)\right)$. Furthermore the closure of all functions constructed in this way is equal to $C\left(\partial\left(\mathbb{S}_{+}^{n-1}\times\mathbb{S}_{+}^{n-1}\right),\mathcal{\mathcal{K}}\left(L^{2}(\mathbb{R}_{+})\right)\right)$
\end{prop}

\begin{proof}
It remains only to prove the statement about the closure of these sets of functions. We will prove only $(ii)$ as $(i)$ and $(iii)$ follow similarly. Let $\varphi\in C_{c}^{\infty}(\mathbb{R}_{+})$ and consider the operator $\left\langle \overline{\varphi}\right|\in\mathcal{B}(L^{2}(\mathbb{R}_{+}),\mathbb{C})$, where we are using the Dirac notation: $\left\langle \overline{\varphi}\right|(u)=\int\varphi(x)u(x)dx$. Let $p\in S_{cl}^{0,0}(\mathbb{R}^{n-1}\times\mathbb{R}^{n-1})$, then we can define a function $t\in S_{cl}^{0,0}(\mathbb{R}^{n-1},\mathbb{R}^{n-1},\mathcal{S}_{+})$ by $t(x',x_{n},\xi')=\left\langle \xi' \right\rangle^{\frac{1}{2}}p(x',\xi')\varphi(\left\langle \xi' \right\rangle x_{n})$. 

This function is such that $t_{.,(0)}(x',\xi',D)\kappa_{\left\langle \xi' \right\rangle}=p_{.,(0)}(x',\xi')\left\langle \overline{\varphi}\right|$ and $t_{(0),.}(x',\xi',D)\kappa_{\left|\xi'\right|}=p_{(0),.}(x',\xi')\left\langle \overline{\varphi}\right|$. The set of these functions is clearly equal to $C^{\infty}\left(\partial\left(\mathbb{S}_{+}^{n-1}\times\mathbb{S}_{+}^{n-1}\right)\right)\otimes C_{c}^{\infty}(\mathbb{R}_{+})$ But this set is dense on $C\left(\partial\left(\mathbb{S}_{+}^{n-1}\times\mathbb{S}_{+}^{n-1}\right),\mathcal{B}(L^{2}(\mathbb{R}_{+}),\mathbb{C})\right)$ - we are making the identification $\mathcal{B}(L^{2}(\mathbb{R}_{+}),\mathbb{C})$ with $L^{2}(\mathbb{R}_{+})$.
In $(iii)$ we just observe that the compact operators $\mathcal{K}(L^{2}(\mathbb{R}_{+}))$
appear because operators with kernel in $\mathcal{S}(\mathbb{R}_{+}\times\mathbb{R}_{+})$
form a dense set of compact operators on $L^{2}(\mathbb{R}_{+})$.
\end{proof}

Also for the pseudodifferential operator there is a similar result.
\begin{prop}
Let $p\in S_{cl}^{0,0}(\mathbb{R}^{n}\times\mathbb{R}^{n})_{tr}$.
Then the function $(x',\xi')\mapsto \kappa_{\left\langle \xi' \right\rangle^{-1}}p(x',0,\xi',D)_{+}\kappa_{\left\langle \xi' \right\rangle}$
belongs to \[
S_{cl}^{0,0}(\mathbb{R}^{n-1},\mathbb{R}^{n-1},\mathcal{B}(L^{2}(\mathbb{R}_{+}))).\]
\end{prop}
\begin{proof}
First we observe that it is an SG-symbol. We know that\[
\left(\kappa_{\left\langle \xi' \right\rangle^{-1}}p(x',0,\xi',D)_{+}\kappa_{\left\langle \xi' \right\rangle}\right)u(x_{n})=p(x',0,\xi',\left\langle \xi' \right\rangle D)_{+}u(x_{n})\]
and we note that\[
\partial_{x'}^{\beta}\partial_{\xi'}^{\alpha}\left(p(x',0,\xi',\left\langle \xi' \right\rangle\xi_{n})\right)=\mbox{sum }a(\xi')\xi_{n}^{k}\left(\partial_{x'}^{\beta}\partial_{\xi'}^{\alpha'}\partial_{\xi_{n}}^{k}p\right)(x',0,\xi',\left\langle \xi' \right\rangle\xi_{n}),\]
where $a\in S^{-l,0}(\mathbb{R}^{n-1}\times\mathbb{R}^{n-1})$ is a symbol independent of $x'$, $|\alpha'|+k+l=|\alpha|$. Finally we note that \foreignlanguage{brazil}{\[
\left\Vert a(\xi')\xi_{n}^{k}\left(\partial_{x'}^{\beta}\partial_{\xi'}^{\alpha'}\partial_{\xi_{n}}^{k}p\right)(x',0,\xi',\left\langle \xi' \right\rangle D)_{+}\right\Vert _{\mathcal{B}(L^{2}(\mathbb{R}_{+}))}=\]
\[
\sup_{\xi_{n}}\left|a(\xi')\xi_{n}^{k}\left(\partial_{x'}^{\beta}\partial_{\xi'}^{\alpha'}\partial_{\xi_{n}}^{k}p\right)(x',0,\xi',\left\langle \xi' \right\rangle \xi_{n})\right|\le C\left\langle \xi' \right\rangle^{-|\alpha|}\left\langle x' \right\rangle^{-|\beta|},\]
}where we used the Toeplitz operator property in the first equality.
That this is true is a consequence of the fact that $p$ satisfies
the transmission property, see, for instance, \cite[Lemma 3.1.5]{Grubbverde}.

The verification that this symbol is classical follows from the observation that it has the required asymptotic expansion in the variables $x'$ and $\xi'$, as the symbols $p$ and $\xi'\mapsto\left\langle \xi'\right\rangle$  are both classical. We can, for instance, compute the first term of the asymptotic expansion in $\xi'$ by \[\left(\kappa_{\left\langle \xi'\right\rangle ^{-1}}p(x',0,\xi',D)_{+}\kappa_{\left\langle \xi'\right\rangle }\right)_{(0),.}=\lim_{\lambda\to\infty}\kappa_{\left\langle \lambda\xi'\right\rangle ^{-1}}p(x',0,\lambda\xi',D)_{+}\kappa_{\left\langle \lambda\xi'\right\rangle }=\] \[\lim_{\lambda\to\infty}p(x',0,\lambda\xi',\left\langle \lambda\xi'\right\rangle D)=p_{(0),.}(x',0,\xi',\left|\xi'\right|D)=\kappa_{\left|\xi'\right|^{-1}}p_{(0),.}(x',0,\xi',D)_{+}\kappa_{\left|\xi'\right|}.\]
and in  $x'$ by
\[\left(\kappa_{\left\langle \xi'\right\rangle ^{-1}}p(x',0,\xi',D)_{+}\kappa_{\left\langle \xi'\right\rangle }\right)_{.,(0)}=\lim_{\lambda\to\infty}\kappa_{\left\langle \xi'\right\rangle ^{-1}}p(\lambda x',0,\xi',D)_{+}\kappa_{\left\langle \xi'\right\rangle }=\kappa_{\left\langle \xi'\right\rangle ^{-1}}p_{.,(0)}(x',0,\xi',D)_{+}\kappa_{\left\langle \xi'\right\rangle }.\]
\end{proof}

\begin{cor}
The functions $\kappa_{|\xi'|^{-1}}p_{(0),.}(x',0,\xi',D)_{+}\kappa_{|\xi'|}$ and $\kappa_{\left\langle \xi' \right\rangle^{-1}}p_{.,(0)}(x',0,\xi',D)_{+}\kappa_{\left\langle \xi' \right\rangle}$ define, in a unique way, a function in $C^{\infty}\left(\partial\left(\mathbb{S}_{+}^{n-1}\times\mathbb{S}_{+}^{n-1}\right),\mathcal{B}\left(L^{2}\left(\mathbb{R}_{+}\right)\right)\right)$.
\end{cor}

Using these results we define the $*$-homomorphism $\gamma$, which we call the boundary principal symbol.

\begin{defn}
The function $\gamma:\mathcal{B}^{(0,0),0}(\mathbb{R}^{n}_{+})\to C^{\infty}\left(\partial\left(\mathbb{S}_{+}^{n-1}\times\mathbb{S}_{+}^{n-1}\right),\mathcal{B}\left(L^{2}\left(\mathbb{R}_{+}\right)\oplus\mathbb{C}\right)\right)$ is the function that assign to each operator $\left(\begin{array}{cc}
op(p)_{+}+op(g) & op(k)\\
op(t) & op(s)\end{array}\right)\in\mathcal{B}^{(0,0),0}(\mathbb{R}^{n}_{+})$ the unique function in $C^{\infty}\left(\partial\left(\mathbb{S}_{+}^{n-1}\times\mathbb{S}_{+}^{n-1}\right),\mathcal{B}\left(L^{2}\left(\mathbb{R}_{+}\right)\oplus\mathbb{C}\right)\right)$ that corresponds - as in Proposition \foreignlanguage{brazil}{\ref{pro:protreze}} and in the previous corollary - to the functions

\[(x',\xi')\mapsto \left(\begin{array}{cc}
\kappa_{|\xi'|^{-1}}p_{(0),.}(x',0,\xi',D)_{+}\kappa_{|\xi'|}+\kappa_{|\xi'|^{-1}}g_{(0),.}(x',\xi',D)\kappa_{|\xi'|} & \kappa_{|\xi'|^{-1}}k_{(0),.}(x',\xi',D)\\
t_{(0),.}(x',\xi',D)\kappa_{|\xi'|} & s_{(0),.}(x',\xi',D)\end{array}\right)\]

and

\[(x',\xi')\mapsto\left(\begin{array}{cc}
\kappa_{\left\langle \xi' \right\rangle^{-1}}p_{.,(0)}(x',0,\xi',D)_{+}\kappa_{\left\langle \xi' \right\rangle}+\kappa_{\left\langle \xi' \right\rangle^{-1}}g_{.,(0)}(x',\xi',D)\kappa_{\left\langle \xi' \right\rangle} & \kappa_{\left\langle \xi' \right\rangle^{-1}}k_{.,(0)}(x',\xi',D)\\
t_{.,(0)}(x',\xi',D)\kappa_{\left\langle \xi' \right\rangle} & s_{.,(0)}(x',\xi',D)\end{array}\right).\]

\end{defn}

We denote by $\left\Vert\gamma(A)\right\Vert$ the supremum over all $\left(z,w\right)\in\partial\left(\mathbb{S}_{+}^{n-1}\times\mathbb{S}_{+}^{n-1}\right)$ of $\left\Vert \gamma(A)(z,w)\right\Vert _{B(L^{2}(\mathbb{R}_{+})\oplus\mathbb{C})}$. We know that 
$\kappa_{\left\langle \xi' \right\rangle}$ and $\kappa_{|\xi'|}$ are unitary. Therefore the matrices below are also unitary: \[
\left(\begin{array}{cc}
\kappa_{\left\langle \xi' \right\rangle} & 0\\
0 & 1\end{array}\right)\mbox{ and }\left(\begin{array}{cc}
\kappa_{|\xi'|} & 0\\
0 & 1\end{array}\right)\]

Hence we conclude that \[
\left\Vert \gamma\left(\left(\begin{array}{cc}
P_{+}+G & K\\
T & S\end{array}\right)\right)\right\Vert =\]

\selectlanguage{brazil}%
\[
\max\left\{ \sup_{(x',\xi')}\left\Vert A_{(0),.}(x',\xi')\right\Vert _{\mathcal{B}(L^{2}(\mathbb{R}_{+})\oplus\mathbb{C})},\sup_{(x',\xi')}\left\Vert A_{.,(0)}(x',\xi')\right\Vert _{\mathcal{B}(L^{2}(\mathbb{R}_{+})\oplus\mathbb{C})}\right\} ,\]

Let us now give another way of understanding the estimate modulo compact of the classical SG-Boutet de Monvel operator.
In order to do that, we need to define some sets.\[
\begin{array}{c}
\mathbb{S}_{++}^{n}=\{z\in\mathbb{R}^{n+1};|z|=1,z_{n}\ge0\,\,\, z_{n+1}\ge0\}.\\
\mathbb{S}_{++}^{n}\times\mathbb{S}_{+}^{n}=\{\left(z,w\right)\in\mathbb{R}^{n+1}\times\mathbb{R}^{n+1};|z|=1,z_{n}\ge0\,\,\, z_{n+1}\ge0\,\,\,\mbox{and}\,\,\,|w|=1,w_{n+1}\ge0\}.\\
\mathbb{S}_{++}^{n}\times\mathbb{S}^{n-1}\cup\mathbb{S}_{+}^{n-1}\times\mathbb{S}_{+}^{n}=\{\left(z,w\right)\in\mathbb{R}^{n+1}\times\mathbb{R}^{n+1};|z|=1,z_{n}\ge0\,\,\, z_{n+1}\ge0\,\,\,\mbox{and}\,\,\,|w|=1,w_{n+1}=0\}\cup\\
\{\left(z,w\right)\in\mathbb{R}^{n+1}\times\mathbb{R}^{n+1};|z|=1,z_{n}\ge0\,\,\, z_{n+1}=0\,\,\,\mbox{and}\,\,\,|w|=1,w_{n+1}\ge0\}.\end{array}\]

We note that using the radial compactification, the interior of the set $\mathbb{S}^{n}_{++}$ is identified with the set $\mathbb{R}^{n}_{+}$. In fact if $x\in\mathbb{R}^{n}\mapsto z=\left(\frac{x}{\left\langle x\right\rangle },\frac{1}{\left\langle x\right\rangle }\right)\in\mathbb{S}^{n}_{+}\subset \mathbb{R}^{n+1}$, then $x_{n}\ge0 \iff z_{n}\ge0$ and $x_{n}=0 \iff z_{n}=0$. We can therefore define the following set:

\begin{defn}
\label{def:sigmapmais} 
The space $C^{\infty}\left(\mathbb{S}_{++}^{n}\times\mathbb{S}^{n-1}\cup\mathbb{S}_{+}^{n-1}\times\mathbb{S}_{+}^{n}\right)_{tr}$ consists of the principal symbols $\sigma(p)$ of functions $p\in S_{cl}^{0,0}(\mathbb{R}^{n}\times\mathbb{R}^{n})_{tr}$, restricted to the space $x_{n}>0$, where we are using the identification of $\mathbb{R}_{+}^{n}$ and $\mathbb{S}_{++}^{n}$, provided by the radial compactification. The restriction of $\sigma(p)$ to $\mathbb{S}_{++}^{n}\times\mathbb{S}^{n-1}\cup\mathbb{S}_{+}^{n-1}\times\mathbb{S}_{+}^{n}$ is denoted by $\sigma(p)^{+}$.
\end{defn}

Using the above definitions, we can write Theorem 12 as: 
 
\begin{cor} \label{cor:teo12}
The estimate modulo compact of the classical SG-Boutet de Monvel operators can be written as \[
\inf_{C\in\mathcal{K}(L^{2}(\mathbb{R}_{+}^{n})\oplus L^{2}(\mathbb{R}^{n-1}))}\left\Vert A+C\right\Vert =\max\left\{\left\Vert \gamma(A)\right\Vert ,\left\Vert \sigma(p)^{+}\right\Vert \right\},\] where $p$ is the symbol of the pseudodifferential part of $A$, and $\left\Vert \sigma(p)^{+}\right\Vert$  is the sup norm on $\mathbb{S}_{++}^{n}\times\mathbb{S}^{n-1}\cup\mathbb{S}_{+}^{n-1}\times\mathbb{S}_{+}^{n}$. \end{cor}

\selectlanguage{english}%
As a consequence of the above estimate, we conclude that $\gamma$ is continuous and can be extended
to the closure of $\mathcal{B}^{(0,0),0}(\mathbb{R}^{n}_{+})$ in $\mathcal{B}\left(L^{2}(\mathbb{R}_{+}^{n})\oplus L^{2}(\mathbb{R}^{n-1})\right)$, which
we denote by $\mathfrak{A}$. We will call this extension $\overline{\gamma}$.
\[
\overline{\gamma}:\mathfrak{A}\to C\left(\partial\left(\mathbb{S}_{+}^{n-1}\times\mathbb{S}_{+}^{n-1}\right),\mathcal{B}(L^{2}(\mathbb{R}_{+})\oplus\mathbb{C})\right).\]

For future reference, we note that $\mathfrak{A}$ can be written as \[
\mathfrak{A}=\left(\begin{array}{cc}
\mathfrak{A}_{11} & \mathfrak{A}_{12}\\
\mathfrak{A}_{21} & \mathfrak{A}_{22}\end{array}\right),\]
where $\mathfrak{A}_{ij}$ is the closure of the entry $(i,j)$ of the matrix in Definition \ref{def:sgboutet}.

Hence, also $\overline{\gamma}$ can be written as a map
of the form \[
\overline{\gamma}\left(\begin{array}{cc}
A_{11} & A_{12}\\
A_{21} & A_{22}\end{array}\right)=\left(\begin{array}{cc}
\overline{\gamma}_{11}\left(A_{11}\right) & \overline{\gamma}_{12}\left(A_{12}\right)\\
\overline{\gamma}_{21}\left(A_{21}\right) & \overline{\gamma}_{22}\left(A_{22}\right)\end{array}\right),\]
where $A_{ij}\in\mathfrak{A}_{ij}$.

\section{The kernel of $\overline{\gamma}$.}

\begin{prop}
\label{pro:kernelgamma}
The kernel of the function $\gamma:\mathcal{B}^{(0,0),0}(\mathbb{R}^{n}_{+})\to C^{\infty}\left(\partial\left(\mathbb{S}_{+}^{n-1}\times\mathbb{S}_{+}^{n-1}\right),\mathcal{B}(L^{2}(\mathbb{R}_{+})\oplus\mathbb{C})\right)$ is equal to the algebra consisting of operators of the form

\begin{equation}
\left(\begin{array}{cc}
P_{+}+G & K\\
T & S\end{array}\right)\label{eq:nome},\end{equation}
 where $G\in\mathcal{G}^{(-1,-1),0}(\mathbb{R}_{+}^{n})$, $K\in\mathcal{K}^{(-1,-1)}(\mathbb{R}_{+}^{n})$,
$T\in\mathcal{T}^{(-1,-1),0}(\mathbb{R}_{+}^{n})$, $S=op(s)$
for some $s\in S_{cl}^{-1,-1}(\mathbb{R}^{n-1}\times\mathbb{R}^{n-1})$ and
$P_{+}=r^{+}op(p)e^{+}$, where $p\in S_{cl}^{0,0}(\mathbb{R}^{n}\times\mathbb{R}^{n})_{tr}$
is such that $p_{(0),.}(x',0,\xi)=0$ and $p_{.,(0)}(x',0,\xi)=0$. That is, the operators $G$, $K$, $T$
and $S$ are of lower order and the principal symbol of $P$ vanishes at the boundary.
\end{prop}
\begin{proof}
If $A$ is an element that satisfies the above conditions, then $A_{(0),.}(x',\xi')=0$ and $A_{.,(0)}(x',\xi')=0$ for all $(x',\xi')$.
Hence $A\in ker(\gamma)$. 

If $\gamma(A)=0$, then the principal symbols
of $t,$ $k$ and $s$ are equal to $0$. Therefore they are all of lower
orders. If for all $(x',\xi')$ the function $p_{(0),.}(x',0,\xi',D)_{+}+g_{(0),.}(x',\xi',D)=0$, then
$p_{(0),.}(x',0,\xi)$ and $g_{(0),.}(x',x_{n},y_{n},\xi')$
are equal to zero, for all $x'$, $\xi'$, $\xi_{n}$, $x_{n}$ and $y_{n}$, as $g_{(0),.}(x',\xi',D)$ is a compact operator
- its kernel belongs to $\mathcal{S}(\mathbb{R}_{++})$ - and $p_{(0),.}(x',0,\xi',D)_{+}$
is a compact Toeplitz operator, hence its symbol vanishes. The same happens with $p_{.,(0)}(x',0,\xi',D)_{_{+}}$
and $g_{.,(0)}(x',\xi',D)$. We conclude that $A$ satisfies the above conditions.
\end{proof}

\begin{defn}

The algebra $\mathcal{J}$ is the closure of $ker\left(\gamma\right)$. It is given by \[
\mathcal{J}=\left(\begin{array}{cc}
\mathcal{J}_{11} & \mathcal{J}_{12}\\
\mathcal{J}_{21} & \mathcal{J}_{22}\end{array}\right),\]
where $\mathcal{J}_{ij}$ is the closure of all operators that appear in the corresponding entry in equation (\ref{eq:nome}).
\end{defn}

We need also to fix some notation: $\left\Vert \overline{\gamma}_{11}(P_{+}+G)\right\Vert $
denotes\foreignlanguage{brazil}{\[
\begin{array}{c}
\max\left\{ \sup_{(x',\xi')}\left\Vert p_{(0),.}(x',0,\xi',D)_{+}+g_{(0),.}(x',\xi',D)\right\Vert _{\mathcal{B}(L^{2}(\mathbb{R}_{+}))},\right.\\
\left.\sup_{(x',\xi')}\left\Vert p_{.,(0)}(x',0,\xi',D)_{+}+g_{.,(0)}(x',\xi',D)\right\Vert _{\mathcal{B}(L^{2}(\mathbb{R}_{+}))}\right\} ,\end{array}\]
}where $\sup_{(x',\xi')}$ means the sup of all $(x',\xi')$ where
the functions are well defined. It is clear that under the usual identification provided by the radial compactification of functions in $\mathbb{R}^{n-1}\times\mathbb{R}^{n-1}$ with functions in $\mathbb{S}_{+}^{n-1}\times\mathbb{S}_{+}^{n-1}$, the above sup can be considered as the supremum over $\partial\left(\mathbb{S}_{+}^{n-1}\times\mathbb{S}_{+}^{n-1}\right)$.

In the same way \foreignlanguage{brazil}{ $\left\Vert \overline{\gamma}_{12}(K)\right\Vert =\sup_{(x',\xi')}\left\Vert \overline{\gamma}_{12}(K)\right\Vert _{\mathcal{B}(\mathbb{C},L^{2}(\mathbb{R}_{+}))}$,
$\left\Vert \overline{\gamma}_{21}(T)\right\Vert =\sup_{(x',\xi')}\left\Vert \overline{\gamma}_{21}(T)\right\Vert _{\mathcal{B}(L^{2}(\mathbb{R}_{+}),\mathbb{C})}$
and $\left\Vert \overline{\gamma}_{22}(S)\right\Vert =\sup_{(x',\xi')}\left|\overline{\gamma}_{22}(S)\right|.$}
\begin{prop}
(Analogous to \cite[Lemma 2]{MeloSchrohe}) Let $P_{+}=op(p)_{+}$,
where $p\in S_{cl}^{0,0}(\mathbb{R}^{n}\times\mathbb{R}^{n})_{tr}$,
and $G\in\mathcal{G}^{(0,0),0}(\mathbb{R}^{n})$. Then there is a constant $c>0$ such that the following
estimate holds\foreignlanguage{brazil}{\[
\inf_{Q_{+}\in\mathcal{Q}}\left\Vert P_{+}+G+Q_{+}\right\Vert _{\mathcal{B}(L^{2}(\mathbb{R}_{+}^{n}))}\le c\left\Vert \overline{\gamma}_{11}(P_{+}+G)\right\Vert ,\]
}where $\mathcal{Q}$ is the set of all operators $Q_{+}=r^{+}op(q)e^{+}$, where 
\textup{$q\in S_{cl}^{0,0}(\mathbb{R}^{n}\times\mathbb{R}^{n})_{tr}$, $q_{(0),.}(x',0,\xi)=0$ and $q_{.,(0)}(x',0,\xi)=0$.} 
\end{prop}
\begin{proof}
Case 1: $P_{+}=0$. Then

\selectlanguage{brazil}%
\[
\inf_{Q_{+}\in\mathcal{Q}}\left\Vert G+Q_{+}\right\Vert _{\mathcal{B}(L^{2}(\mathbb{R}_{+}^{n}))}\le\inf_{C\in\mathcal{K}(L^{2}(\mathbb{R}_{+}^{n}))}\left\Vert G+C\right\Vert _{\mathcal{B}(L^{2}(\mathbb{R}_{+}^{n}))}=\left\Vert \overline{\gamma}_{11}(G)\right\Vert ,\]
\foreignlanguage{english}{because the operators $op(q)_{+}$, with
$q\in S^{-\infty,-\infty}(\mathbb{R}^{n}\times\mathbb{R}^{n})$,
belong to $\mathcal{Q}$ and they form a dense subset of the set of compact operators. The second equality follows from Corollary \ref{cor:teo12}.}

\selectlanguage{english}%

$\vspace{0.2cm}$

Case 2: $G=0$.

Let $\chi \in C^{\infty}_{c}(\mathbb{R})$ be a function such that $0\le\chi\le 1$ and that it is equal to $1$ in a neighborhood of $0$. Let $\tilde{p}\in S_{cl}^{0,0}(\mathbb{R}^{n}\times\mathbb{R}^{n})$ be given by $\tilde{p}(x,\xi)=p(x',0,\xi)\chi\left(\frac{x_{n}}{\left\langle x'\right\rangle}\right)$. If $\tilde{P}=op(\tilde{p})$, then
 $\left(P-\tilde{P}\right)_{+}\in\mathcal{Q}$, as $p_{(0),.}(x',0,\xi)=\tilde{p}_{(0),.}(x',0,\xi)$
and $p_{.,(0)}(x',0,\xi)=\tilde{p}_{.,(0)}(x',0,\xi)$. Therefore\foreignlanguage{brazil}{\[
\inf_{Q_{+}\in\mathcal{Q}}\left\Vert P_{+}+Q_{+}\right\Vert _{\mathcal{B}(L^{2}(\mathbb{R}_{+}^{n}))}=\inf_{Q_{+}\in\mathcal{Q}}\left\Vert P_{+}+(\tilde{P}_{+}-P_{+})+Q_{+}\right\Vert _{\mathcal{B}(L^{2}(\mathbb{R}_{+}^{n}))}=\inf_{Q_{+}\in\mathcal{Q}}\left\Vert \tilde{P}_{+}+Q_{+}\right\Vert _{\mathcal{B}(L^{2}(\mathbb{R}_{+}^{n}))}.\]
} Again, as $\mathcal{Q}$ contains a dense set of compact operators:
\foreignlanguage{brazil}{\[
\inf_{Q_{+}\in\mathcal{Q}}\left\Vert \tilde{P}_{+}+Q_{+}\right\Vert _{\mathcal{B}(L^{2}(\mathbb{R}_{+}^{n}))}\le\inf_{C\in\mathcal{K}(L^{2}(\mathbb{R}^{n}))}\left\Vert \tilde{P}_{+}+C_{+}\right\Vert _{\mathcal{B}(L^{2}(\mathbb{R}_{+}^{n}))}\le\inf_{C\in\mathcal{K}(L^{2}(\mathbb{R}^{n}))}\left\Vert \tilde{P}+C\right\Vert _{\mathcal{B}(L^{2}(\mathbb{R}^{n}))}.\]
}Using Theorem 3, we conclude that\foreignlanguage{brazil}{\[
\inf_{C\in\mathcal{K}(L^{2}(\mathbb{R}^{n}))}\left\Vert \tilde{P}+C\right\Vert _{\mathcal{B}(L^{2}(\mathbb{R}^{n}))}=\max\left\{ \sup_{(x,\xi)}\left|p_{(0),.}(x',0,\xi)\chi\left(\frac{x_{n}}{\left\langle x'\right\rangle}\right)\right|,\sup_{(x,\xi)}\left|p_{.,(0)}(x',0,\xi)\chi\left(\frac{x_{n}}{\left|x'\right|}\right)\right|\right\} \le\]
\[
\max\left\{ \sup_{(x',\xi')}\sup_{\xi_{n}}\left|p_{(0),.}(x',0,\xi)\right|,\sup_{(x',\xi')}\sup_{\xi_{n}}\left|p_{.,(0)}(x',0,\xi)\right|\right\} =\]
\[
\max\left\{ \sup_{(x',\xi')}\left\Vert p_{(0),.}(x',0,\xi',D)_{+}\right\Vert _{\mathcal{B}(L^{2}(\mathbb{R}_{+}))},\sup_{(x',\xi')}\left\Vert p_{.,(0)}(x',0,\xi',D)_{+}\right\Vert _{\mathcal{B}(L^{2}(\mathbb{R}_{+}))}\right\} =\left\Vert \overline{\gamma}_{11}(P_{+})\right\Vert.\]
}In the second equality, we used the Toeplitz operators properties \cite[Lemma 3.1.5]{Grubbverde}.

$\vspace{0.2cm}$

Case 3: General case.

Again by \cite[Lemma 3.1.5]{Grubbverde} and because $g_{(0),.}(x',\xi',D)$ is compact, we have \foreignlanguage{brazil}{\[
\left\Vert \gamma_{11}(P_{+})\right\Vert =\max\left\{ \sup_{(x',\xi')}\left\Vert p_{(0),.}(x',0,\xi',D)_{+}\right\Vert _{\mathcal{B}(L^{2}(\mathbb{R}_{+}))},\sup_{(x',\xi')}\left\Vert p_{.,(0)}(x',0,\xi',D)_{+}\right\Vert _{\mathcal{B}(L^{2}(\mathbb{R}_{+}))}\right\} =\]
\[
\max\left\{ \sup_{(x',\xi')}\left(\inf_{C\in\mathcal{K}(L^{2}(\mathbb{R}_{+}))}\left\Vert p_{(0),.}(x',0,\xi',D)_{+}+C\right\Vert _{\mathcal{B}(L^{2}(\mathbb{R}_{+}))}\right),\right.\]
\[
\left.\sup_{(x',\xi')}\left(\inf_{C\in\mathcal{K}(L^{2}(\mathbb{R}_{+}))}\left\Vert p_{.,(0)}(x',0,\xi',D)_{+}+C\right\Vert _{\mathcal{B}(L^{2}(\mathbb{R}_{+}))}\right)\right\} \le\]
\[
\max\left\{ \sup_{(x',\xi')}\left\Vert p_{(0),.}(x',0,\xi',D)_{+}+g_{(0),.}(x',\xi',D)\right\Vert _{\mathcal{B}(L^{2}(\mathbb{R}_{+}))},\right.\]
\[
\left.\sup_{(x',\xi')}\left\Vert p_{.,(0)}(x',0,\xi',D)_{+}+g_{.,(0)}(x',\xi',D)\right\Vert _{\mathcal{B}(L^{2}(\mathbb{R}_{+}))}\right\} =\left\Vert \overline{\gamma}_{11}(P_{+}+G)\right\Vert .\]
}

Hence \foreignlanguage{brazil}{\[
\left\Vert \overline{\gamma}_{11}(G)\right\Vert  \le\left\Vert \overline{\gamma}_{11}(P_{+}+G)\right\Vert +\left\Vert \overline{\gamma}_{11}(P_{+})\right\Vert \le2\left\Vert \overline{\gamma}_{11}(P_{+}+G)\right\Vert .\]
}

Finally we conclude that\foreignlanguage{brazil}{\[
\inf_{Q_{+}\in\mathcal{Q}}\left\Vert P_{+}+G+Q_{+}\right\Vert _{\mathcal{B}(L^{2}(\mathbb{R}_{+}^{n}))}\le\]
\[
\inf_{Q_{1+}\in\mathcal{Q}}\left\Vert P_{+}+Q_{1+}\right\Vert _{\mathcal{B}(L^{2}(\mathbb{R}_{+}^{n}))}+\inf_{Q_{2+}\in\mathcal{Q}}\left\Vert G+Q_{2+}\right\Vert _{\mathcal{B}(L^{2}(\mathbb{R}_{+}^{n}))}\le\left\Vert \overline{\gamma}_{11}(P_{+})\right\Vert +\left\Vert \overline{\gamma}_{11}(G)\right\Vert \le3\left\Vert \overline{\gamma}_{11}(P_{+}+G)\right\Vert .\]
}
\end{proof}

\begin{cor}
The following estimate holds for any $A\in\mathfrak{A}_{11}$:\foreignlanguage{brazil}{\[
\inf_{A'\in\mathcal{J}_{11}}\left\Vert A+A'\right\Vert _{\mathcal{B}(L^{2}(\mathbb{R}_{+}^{n}))}\le c\left\Vert \overline{\gamma_{11}}(A)\right\Vert .\]
}\end{cor}
\begin{proof}
We remark that $\mathcal{Q}\subset\mathcal{J}_{11}$. Therefore for
any $P_{+}+G$ we have that \foreignlanguage{brazil}{\[
\inf_{A'\in\mathcal{J}_{11}}\left\Vert P_{+}+G+A'\right\Vert _{\mathcal{B}(L^{2}(\mathbb{R}_{+}^{n}))}\le\inf_{Q_{+}\in\mathcal{Q}}\left\Vert P_{+}+G+Q_{+}\right\Vert _{\mathcal{B}(L^{2}(\mathbb{R}_{+}^{n}))}\le c\left\Vert \overline{\gamma_{11}}(P_{+}+G)\right\Vert .\]
}

Now we only have to use that $\overline{\gamma}_{11}:\mathfrak{A}_{11}\to C\left(\partial\left(\mathbb{S}_{+}^{n-1}\times\mathbb{S}_{+}^{n-1}\right),\mathcal{B}(L^{2}(\mathbb{R}_{+}))\right)$
is a continuous map and that the operators $P_{+}+G$ form a dense subset of
$\mathfrak{A}_{11}$.
\end{proof}
Finally the main result of all this discussion is:
\begin{cor}
$ker\left(\overline{\gamma}\right)=\mathcal{J}:=\overline{ker\left(\gamma\right)}$.\end{cor}
\begin{proof}
It is obvious that $ker\left(\gamma\right)\subset ker\left(\overline{\gamma}\right)$.
Hence $\mathcal{J}\subset ker\left(\overline{\gamma}\right)$, as $ker\left(\overline{\gamma}\right)$
is closed. In order to prove that $ker\left(\overline{\gamma}\right)\subset\mathcal{J}$,
let us denote $A\in ker\left(\overline{\gamma}\right)$ by\[
A=\left(\begin{array}{cc}
A_{11} & A_{12}\\
A_{21} & A_{22}\end{array}\right),\]

We know that the regularizing
operators belong to $\mathcal{J}$ and, therefore, also the compact operators. By Corollary \ref{cor:teo12}, we conclude that

\selectlanguage{brazil}%
\[
\inf_{A_{12}'\in\mathcal{J}_{12}}\left\Vert A_{12}+A_{12}'\right\Vert _{\mathcal{B}(L^{2}(\mathbb{R}^{n-1}),L^{2}(\mathbb{R}_{+}^{n}))}\le\inf_{C\in\mathcal{K}(L^{2}(\mathbb{R}^{n-1}),L^{2}(\mathbb{R}_{+}^{n}))}\left\Vert A_{12}+C\right\Vert _{\mathcal{B}(L^{2}(\mathbb{R}^{n-1}),L^{2}(\mathbb{R}_{+}^{n}))}\le\left\Vert \overline{\gamma}_{12}(A_{12})\right\Vert ,\]

\[
\inf_{A_{21}'\in\mathcal{J}_{21}}\left\Vert A_{21}+A_{21}'\right\Vert _{\mathcal{B}(L^{2}(\mathbb{R}_{+}^{n}),L^{2}(\mathbb{R}^{n-1}))}\le\inf_{C\in\mathcal{K}(L^{2}(\mathbb{R}_{+}^{n}),L^{2}(\mathbb{R}^{n-1}))}\left\Vert A_{21}+C\right\Vert _{\mathcal{B}(L^{2}(\mathbb{R}_{+}^{n}),L^{2}(\mathbb{R}^{n-1}))}\le\left\Vert \overline{\gamma}_{21}(A_{21})\right\Vert ,\]

\[
\inf_{A_{22}'\in\mathcal{J}_{22}}\left\Vert A_{22}+A_{22}'\right\Vert _{\mathcal{B}(L^{2}(\mathbb{R}^{n-1}),L^{2}(\mathbb{R}^{n-1}))}\le\inf_{C\in\mathcal{K}(L^{2}(\mathbb{R}^{n-1}),L^{2}(\mathbb{R}^{n-1}))}\left\Vert A_{22}+C\right\Vert _{\mathcal{B}(L^{2}(\mathbb{R}^{n-1}),L^{2}(\mathbb{R}^{n-1}))}\le\left\Vert \overline{\gamma}_{22}(A_{22})\right\Vert .\]

\selectlanguage{english}%
As $\overline{\gamma}_{jk}(A_{jk})=0$,
we conclude that $A_{jk}\in\mathcal{J}_{jk}$ if $(j,k)\ne(1,1)$. For
$(j,k)=(1,1)$ we have\[
\inf_{A_{11}'\in\mathcal{J}_{11}}\left\Vert A_{11}+A'_{11}\right\Vert _{\mathcal{B}(L^{2}(\mathbb{R}_{+}^{n}))}\le c\left\Vert \overline{\gamma}_{11}(A_{11})\right\Vert .\]

As $\overline{\gamma}_{11}(A_{11})=0$, we finally conclude that $A\in\mathcal{J}$.
\end{proof}

\section{The image of $\overline{\gamma}$.}\label{sec:image}

\begin{defn}
We define the set $\mathfrak{T}$ as the $C^{*}$ algebra of bounded
operators on $L^{2}(\mathbb{R}_{+})$ generated by $\{p(D)_{+};p\in\mathcal{H}_{0}\}$,
where $\mathcal{H}_{0}$ was defined on page \pageref{H}. Let
$u\in L^{2}(\mathbb{R}_{+})$, then $p(D)_{+}$ is defined by  \[
p(D)_{+}u=r^{+}\left(\frac{1}{2\pi}\int e^{-it\xi}p(\xi)\widehat{e^{+}u}(\xi)d\xi\right).\]

\end{defn}
Using that \cite[Lemma 3.1.5]{Grubbverde}
\[\inf_{C\in\mathcal{K}(L^{2}(\mathbb{R}_{+}))}\left\Vert p(D)_{+}+C\right\Vert =\sup_{\xi_{n}\in\mathbb{R}}\left|p(\xi_{n})\right|\] and the fact that $\mathfrak{T}$ has compact commutators, we see that the assignment $p\mapsto p(\infty):=\lim_{\left|\xi_{n}\right|\to\infty}p(\xi_{n})\in\mathbb{C}$ extends to a $C^{*}$-algebra homomorphism, denoted by $\lambda:\mathfrak{T}\to\mathbb{C}$.

\begin{defn}
We define the $C^{*}$-algebra $\mathfrak{T}_{0}$ as the kernel of $\lambda:\mathfrak{T}\to\mathbb{C}$.
\end{defn}
The algebra $\mathfrak{T}$ is unitarily equivalent to the algebra
of Toeplitz operators of continuous symbols as observed in \cite{MeloSchrohe}.
Both algebras $\mathfrak{T}$ and $\mathfrak{T}_{0}$ contain the
compact operators of $L^{2}(\mathbb{R}_{+})$.

We have already seen that\[
\begin{array}{c}
\mbox{Im}\overline{\gamma}_{12}=C\left(\partial\left(\mathbb{S}_{+}^{n-1}\times\mathbb{S}_{+}^{n-1}\right),\mathcal{B}(\mathbb{C},L^{2}(\mathbb{R}_{+}))\right)\\
\mbox{Im}\overline{\gamma}_{21}=C\left(\partial\left(\mathbb{S}_{+}^{n-1}\times\mathbb{S}_{+}^{n-1}\right),\mathcal{B}(L^{2}(\mathbb{R}_{+}),\mathbb{C})\right)\\
\mbox{Im}\overline{\gamma}_{22}=C\left(\partial\left(\mathbb{S}_{+}^{n-1}\times\mathbb{S}_{+}^{n-1}\right)\right).\end{array}\]

The only component that is still not clear is $\mbox{Im}(\overline{\gamma}_{11})$.
\begin{thm}
\label{thm:imagem de gamma}As a Banach space $\mbox{Im}(\overline{\gamma}_{11})$ is
isomorphic to \textup{$C\left(\partial\left(\mathbb{S}_{+}^{n-1}\times\mathbb{S}_{+}^{n-1}\right),\mathfrak{T}_{0} \right)\oplus C\left(S_{+x'}^{n-1}\right)$.} 
\end{thm}
In the above theorem and in what follows, we denote and identify the set
$C(\mathbb{S}_{+x'}^{n-1})$ with the subset of functions in $C\left(\partial\left(\mathbb{S}_{+}^{n-1}\times\mathbb{S}_{+}^{n-1}\right)\right)$
that only depends on the first copy of $\mathbb{S}^{n-1}_{+}$.

The sets $C(\mathbb{S}_{+x'}^{n-1})$ and $C\left(\partial\left(\mathbb{S}_{+}^{n-1}\times\mathbb{S}_{+}^{n-1}\right)\right)$
will also denote subsets of $C\left(\partial\left(\mathbb{S}_{+}^{n-1}\times\mathbb{S}_{+}^{n-1}\right),\mathfrak{T}\right)$
which are multiplications of functions in $C(\mathbb{S}_{+x'}^{n-1})$
and $C\left(\partial\left(\mathbb{S}_{+}^{n-1}\times\mathbb{S}_{+}^{n-1}\right)\right)$,
respectively, by the identity of $\mathcal{B}(L^{2}(\mathbb{R}_{+}))$.

As $\overline{\gamma}_{11}$ is a $C^{*}$-algebra homomorphism, the set $\mbox{Im}(\overline{\gamma}_{11})$ is closed. Hence, in order
to describe $\mbox{Im}(\overline{\gamma}_{11})$, it is enough to
find the closure of the set that consists of elements of the form
$\overline{\gamma}_{11}(P_{+}+G)$, which are identified with the following pairs of functions:
\[
\left(\kappa_{|\xi'|^{-1}}\left(p_{(0),.}(x',0,\xi',D)_{+}+g_{(0),.}(x',\xi',D)\right)\kappa_{|\xi'|},\kappa_{\left\langle \xi' \right\rangle^{-1}}\left(p_{.,(0)}(x',0,\xi',D)_{+}+g_{.,(0)}(x',\xi',D)\right)\kappa_{\left\langle \xi' \right\rangle}\right),\]
where $p\in S_{cl}^{0,0}(\mathbb{R}^{n}\times\mathbb{R}^{n})_{tr}$
and $g\in S_{cl}^{0,0}(\mathbb{R}^{n-1},\mathbb{R}^{n-1},\mathcal{S}_{++})$.

Let us prove the Theorem step by step.
\begin{lem}
\label{lem:o fecho contem C(S,T0)}The closure of the set \[
\left(\kappa_{|\xi'|^{-1}}p_{(0),.}(x',0,\xi',D)_{+}\kappa_{|\xi'|},\kappa_{\left\langle \xi' \right\rangle^{-1}}p_{.,(0)}(x',0,\xi',D)_{+}\kappa_{\left\langle \xi' \right\rangle}\right).\]
 contains $C\left(\partial\left(\mathbb{S}_{+}^{n-1}\times\mathbb{S}_{+}^{n-1}\right),\mathfrak{T}_{0}\right)$.\end{lem}
\begin{proof}
Let $\tilde{p}\in S_{cl}^{0,0}(\mathbb{R}^{n-1}\times\mathbb{R}^{n-1})$. Let us choose functions $\varphi\in C_{c}^{\infty}(\mathbb{R})$ and $\chi\in C_{c}^{\infty}(\mathbb{R})$
such that $\chi=1$ and $\varphi=1$ in a neighborhood of 0. We define
$p\in S_{cl}^{0,0}(\mathbb{R}^{n}\times\mathbb{R}^{n})$ by $p(x,\xi):=\tilde{p}(x',\xi')\chi\left(\frac{x_{n}}{\left\langle x'\right\rangle}\right)\varphi\left(\frac{\xi_{n}}{\left\langle \xi' \right\rangle}\right)$.
This symbol satisfies the transmission property and it is such that
\[
\left(\kappa_{|\xi'|^{-1}}p_{(0),.}(x',0,\xi',D)_{+}\kappa_{|\xi'|},\kappa_{\left\langle \xi' \right\rangle^{-1}}p_{.,(0)}(x',0,\xi',D)_{+}\kappa_{\left\langle \xi' \right\rangle}\right)=\]
\[
\left(\kappa_{|\xi'|^{-1}}\tilde{p}_{(0),.}(x',\xi')\varphi\left(\frac{D}{|\xi'|}\right)_{+}\kappa_{|\xi'|},\kappa_{\left\langle \xi' \right\rangle^{-1}}\tilde{p}_{.,(0)}(x',\xi')\varphi\left(\frac{D}{\left\langle \xi' \right\rangle}\right)_{+}\kappa_{\left\langle \xi' \right\rangle}\right)=\]
\[
\left(\tilde{p}_{(0),.}(x',\xi')\varphi\left(D\right)_{+},\tilde{p}_{.,(0)}(x',\xi')\varphi\left(D\right)_{+}\right).\]
Using the notation $C_{c}^{\infty}(\mathbb{R})(D)=\{\varphi(D)_{+};\varphi\in C_{c}^{\infty}(\mathbb{R})\}$,
we conclude that $C^{\infty}\left(\partial\left(\mathbb{S}_{+}^{n-1}\times\mathbb{S}_{+}^{n-1}\right)\right)\otimes C_{c}^{\infty}(\mathbb{R})(D)$
belongs to the image of $\overline{\gamma}_{11}$.
As $C^{\infty}\left(\partial\left(\mathbb{S}_{+}^{n-1}\times\mathbb{S}_{+}^{n-1}\right)\right)$
is dense in $C\left(\partial\left(\mathbb{S}_{+}^{n-1}\times\mathbb{S}_{+}^{n-1}\right)\right)$
and $C_{c}^{\infty}(\mathbb{R})(D)$ is dense in $\mathfrak{T}_{0}$, the above functions form a dense set of
\[
C\left(\partial\left(\mathbb{S}_{+}^{n-1}\times\mathbb{S}_{+}^{n-1}\right)\right)\hat{\otimes}\mathfrak{T}_{0}=C\left(\partial\left(\mathbb{S}_{+}^{n-1}\times\mathbb{S}_{+}^{n-1}\right),\mathfrak{T}_{0}\right).\]
\end{proof}

\begin{lem}
$\mbox{Im}(\overline{\gamma}_{11})\cap C\left(\partial\left(\mathbb{S}_{+}^{n-1}\times\mathbb{S}_{+}^{n-1}\right)\right)=C(\mathbb{S}_{+x'}^{n-1})$.\end{lem}
\begin{proof}
Let us first prove $\supset$. We fix a function $\varphi \in C^{\infty}_{c}(\mathbb{R})$ that is equal to $1$ in a neighborhood of $0$ and choose $\tilde{p}\in S_{cl}^{0,0}(\mathbb{R}^{n-1}\times\mathbb{R}^{n-1})$ that does not depend on $\xi'$. We define $p\in S^{0,0}_{cl}(\mathbb{R}^{n}\times\mathbb{R}^{n})$ by
$p(x):=\tilde{p}(x')\varphi\left(\frac{x_{n}}{\left\langle x'\right\rangle}\right)$. Therefore
 $p_{(0),.}(x',0)=\tilde{p}(x')$ and $p_{.,(0)}(x',0)=\tilde{p}_{.,(0)}(x')$. The function $\overline{\gamma_{11}}\left(op(p)_{+}\right)$ is identified with the pair of functions of $(x',\xi')$: \[
\left(\kappa_{|\xi'|^{-1}}\left(p_{(0),.}(x',0,\xi',D)_{+}\right)\kappa_{|\xi'|},\kappa_{\left\langle \xi' \right\rangle^{-1}}\left(p_{.,(0)}(x',0,\xi',D)_{+}\right)\kappa_{\left\langle \xi' \right\rangle}\right)=\]
\[
\left(\tilde{p}(x'),\tilde{p}_{.,(0)}(x')\right),\]
we then conclude that $C^{\infty}(\mathbb{S}_{+x'}^{n-1})$, and hence $C(\mathbb{S}_{+x'}^{n-1})$, is contained in $\mbox{Im}\left(\overline{\gamma}_{11}\right)\cap C\left(\partial\left(\mathbb{S}_{+}^{n-1}\times\mathbb{S}_{+}^{n-1}\right)\right)$.

Now let us prove $\subset$. Suppose that $f\in\mbox{Im}(\overline{\gamma}_{11})\cap C\left(\partial\left(\mathbb{S}_{+}^{n-1}\times\mathbb{S}_{+}^{n-1}\right)\right)$.
Then there are functions $\tilde{f}_{(0),.}\in C(\mathbb{R}^{n-1}\times\mathbb{S}^{n-2})$
and $\tilde{f}_{.,(0)}\in C(\mathbb{S}^{n-2}\times\mathbb{R}^{n-1})$
that corresponds to $f$, that is \[
\tilde{f}_{(0),.}=f\circ(RC\times i)\]
 and \[
\tilde{f}_{.,(0)}=f\circ(i\times RC),\]
where $i:\mathbb{S}^{n-2}=\left\{ z\in\mathbb{R}^{n-1},\,\left|z\right|=1\right\} \to \mathbb{S}^{n-2}=\left\{ z\in\mathbb{R}^{n},\,\left|z\right|=1 , z_{n}=0\right\}$ is the "identity". As $f\in\mbox{Im}(\overline{\gamma})$, we
know that for any $\delta>0$, there exists $p\in S_{cl}^{0,0}(\mathbb{R}^{n}\times\mathbb{R}^{n})_{tr}$
that satisfies the transmission property and $g\in S_{cl}^{0,0}(\mathbb{R}_{x'}^{n-1},\mathbb{R}_{\xi'}^{n-1},\mathcal{S}_{++})$
such that\[
\begin{array}{c}
\left\Vert \kappa_{|\xi'|^{-1}}g_{(0),.}(x',\xi',D)\kappa_{|\xi'|}+\kappa_{|\xi'|^{-1}}p_{(0),.}(x',0,\xi',D)_{+}\kappa_{|\xi'|}-\tilde{f}_{(0),.}\left(x',\frac{\xi'}{|\xi'|}\right)I\right\Vert _{\mathcal{B}(L^{2}(\mathbb{R}_{+}))}<\delta.\\
\left\Vert \kappa_{\left\langle \xi' \right\rangle^{-1}}g_{.,(0)}(x',\xi',D)\kappa_{\left\langle \xi' \right\rangle}+\kappa_{\left\langle \xi' \right\rangle^{-1}}p_{.,(0)}(x',0,\xi',D)_{+}\kappa_{\left\langle \xi' \right\rangle}-\tilde{f}_{.,(0)}\left(\frac{x'}{|x'|},\xi'\right)I\right\Vert _{\mathcal{B}(L^{2}(\mathbb{R}_{+}))}<\delta.\end{array}\]

The operators $g_{(0),.}(x',\xi',D)$ and $g_{.,(0)}(x',\xi',D)$ are compact. Hence
\[
\inf_{C\in\mathcal{K}(L^{2}(\mathbb{R}_{+}))}\left\Vert \kappa_{|\xi'|^{-1}}p_{(0),.}(x',0,\xi',D)_{+}\kappa_{|\xi'|}-\tilde{f}_{(0),.}\left(x',\frac{\xi'}{|\xi'|}\right)I+C\right\Vert _{\mathcal{B}(L^{2}(\mathbb{R}_{+}))}<\delta\]
\[
\Longrightarrow\inf_{C\in\mathcal{K}(L^{2}(\mathbb{R}_{+}))}\left\Vert p_{(0),.}(x',0,\xi',D)_{+}-\tilde{f}_{(0),.}\left(x',\frac{\xi'}{|\xi'|}\right)I+C\right\Vert _{\mathcal{B}(L^{2}(\mathbb{R}_{+}))}<\delta,\]
where we used that $\kappa_{|\xi'|}$ and $\kappa_{|\xi'|^{-1}}$
are unitary.

Using the Toeplitz property, we conclude that $\sup_{\xi_{n}\in\mathbb{R}}\left|p_{(0),.}(x',0,\xi',\xi_{n})-\tilde{f}_{(0),.}\left(x',\frac{\xi'}{|\xi'|}\right)\right|<\delta$.
This implies that \[
\lim_{\xi_{n}\to\infty}\left|p_{(0),.}(x',0,\xi',\xi_{n})-\tilde{f}_{(0),.}\left(x',\frac{\xi'}{|\xi'|}\right)\right|=\left|p_{(0),.}(x',0,0,1)-\tilde{f}_{(0),.}\left(x',\frac{\xi'}{|\xi'|}\right)\right|<\delta\]
As $\left((x',\xi')\mapsto p_{(0),.}(x',0,0,1)\right)\in S_{cl}^{0,0}(\mathbb{R}^{n-1}\times \mathbb{R}^{n-1})$ does not depend on $\xi'$, we
conclude that for any $\delta>0$, there exists a $g\in C^{\infty}(\mathbb{S}_{+x'}^{n-1})$,
such that $g\circ RC=p_{(0),.}(x',0,0,1)$ and $\left|\tilde{f}_{(0),.}\left(x',\frac{\xi'}{\left|\xi'\right|}\right)-g(RC(x'))\right|<\delta$, for all $(x',\xi')$.
Hence $\tilde{f}_{(0),.}$ does not depend on $\xi'$.

We can use the same argument for $\tilde{f}_{.,(0)}$ and conclude that it also does not depend on $\xi'$.

\end{proof}
We are finally in position to prove Theorem \ref{thm:imagem de gamma}.
\begin{proof}
(Theorem \ref{thm:imagem de gamma}). Let us define $\tilde{\lambda}:C\left(\partial\left(\mathbb{S}_{+}^{n-1}\times\mathbb{S}_{+}^{n-1}\right),\mathfrak{T}\right)\to C\left(\partial\left(\mathbb{S}_{+}^{n-1}\times\mathbb{S}_{+}^{n-1}\right)\right)$
by the following formula\[
\left(\tilde{\lambda}f\right)(z,w)=\lambda\left(f(z,w)\right).\]

Let $H\in\mbox{Im}(\overline{\gamma}_{11})$. Then $H\in C\left(\partial\left(\mathbb{S}_{+}^{n-1}\times\mathbb{S}_{+}^{n-1}\right),\mathfrak{T}\right)$.
This means that \[
H-\tilde{\lambda}(H)\in C\left(\partial\left(\mathbb{S}_{+}^{n-1}\times\mathbb{S}_{+}^{n-1}\right),\mathfrak{T}_{0}\right)\]
By Lemma \ref{lem:o fecho contem C(S,T0)}, we conclude that $H-\tilde{\lambda}(H)\in\mbox{Im}(\overline{\gamma}_{11})$. Hence \[
\tilde{\lambda}(H)=H-(H-\tilde{\lambda}(H))\in C\left(\partial\left(\mathbb{S}_{+}^{n-1}\times\mathbb{S}_{+}^{n-1}\right)\right)\cap\mbox{Im}(\overline{\gamma}_{11})=C(\mathbb{S}_{+x'}^{n-1}).\]
As
$H=\left(H-\tilde{\lambda}(H)\right)+\tilde{\lambda}(H)$, we conclude
that\[
\mbox{Im}(\overline{\gamma}_{11})=C\left(\partial\left(\mathbb{S}_{+}^{n-1}\times\mathbb{S}_{+}^{n-1}\right),\mathfrak{T}_{0}\right)+C(\mathbb{S}_{+x'}^{n-1}).\]
As $\lambda$ vanishes on $C\left(\partial\left(\mathbb{S}_{+}^{n-1}\times\mathbb{S}_{+}^{n-1}\right),\mathfrak{T}_{0}\right)$ and only vanishes on a $g\in C\left(\mathbb{S}_{+}^{n-1}\right)$ if $g\equiv0$, we conclude that \[
C\left(\partial\left(\mathbb{S}_{+}^{n-1}\times\mathbb{S}_{+}^{n-1}\right),\mathfrak{T}_{0}\right)\cap C(S_{+x'}^{n-1})=\{0\}.\]

\end{proof}

\section{The K-Theory of the algebra.}

Since the set of all regularizing operators (integral operators with smooth kernels) is dense in the ideal $\mathcal{K}$ of all compact operators of $\mathcal{B}(L^{2}(\mathbb{R}_{+}^{n})\oplus L^{2}(\mathbb{R}^{n-1}))$ and, by Proposition \ref{pro:kernelgamma}, these regularizing operators belong to $\mathcal{J}=\ker(\overline{\gamma})$, we conclude that $\mathcal{K}\subset\mathcal{J}$.

For the computation of the K-groups of the algebra $\mathfrak{A}/\mathcal{K}$,
where $\mathfrak{A}$ is the closure of the algebra $\mathcal{B}^{(0,0),0}(\mathbb{R}^{n}_{+})$
in $\mathcal{B}(L^{2}(\mathbb{R}_{+}^{n})\oplus L^{2}(\mathbb{R}^{n-1}))$, we are going to use the following exact sequence induced by $\overline{\gamma}$:

\begin{equation}
0\to\frac{\mathcal{J}}{\mathcal{K}}\to\frac{\mathfrak{A}}{\mathcal{K}}\to\frac{\mathfrak{A}}{\mathcal{J}}\to0.
\label{eq:basica}
\end{equation}

\selectlanguage{english}%

\selectlanguage{brazil}%
\begin{prop}
\label{pro:kergamma em simb prin de p}

\selectlanguage{english}%

There is an injective $C^{*}$-algebra homomorphism \foreignlanguage{brazil}{\[
j:\frac{\mathcal{J}}{\mathcal{K}}\to C(\mathbb{S}_{++}^{n}\times\mathbb{S}^{n-1}\cup\mathbb{S}_{+}^{n-1}\times\mathbb{S}_{+}^{n}),\]
}such that for $A=\left(\begin{array}{cc}
P_{+}+G & K\\
T & S\end{array}\right)$, $j(A+\mathcal{K})=\sigma(p)^{+}$, where $P=op(p)$.
Moreover, let $z\in\mathbb{R}^{n+1}$ be written as $(z',z_{n},z_{n+1})\in\mathbb{R}^{n-1}\times\mathbb{R}\times\mathbb{R}$ and let $w$ belong to $\mathbb{R}^{n+1}$. Hence \[\mbox{Im}(j)=\left\{ \varphi\in C\left(\mathbb{S}_{++}^{n}\times\mathbb{S}^{n-1}\cup\mathbb{S}_{+}^{n-1}\times\mathbb{S}_{+}^{n}\right);\varphi(z',0,z_{n+1},w)=0\right\}\] \end{prop}
\selectlanguage{english}%

Remark:
The set $\mathbb{S}_{++}^{n}\times\mathbb{S}^{n-1}$ above can be regarded as the co-sphere bundle (with spheres of infinite radius) of a compactification of $\mathbb{R}^{n}_{+}$, while the set $\mathbb{S}_{+}^{n-1}\times\mathbb{S}_{+}^{n}$ corresponds to the points over $|x|=\infty$ that are needed in the picture to take care of the behavior of the symbols for large $|x|$. Their intersection is the set where both $|x|$ and $|\xi|$ are infinite.

\begin{proof}
Let $f:\mbox{ker}(\gamma)\to C^{\infty}\left(\mathbb{S}_{++}^{n}\times\mathbb{S}^{n-1}\cup\mathbb{S}_{+}^{n-1}\times\mathbb{S}_{+}^{n}\right)_{tr}$
be given by \[
f\left(\begin{array}{cc}
P_{+}+G & K\\
T & S\end{array}\right)=\sigma(p)^{+},\]
where $P_{+}=r^{+}op(p)e^{+}$ and $\sigma(p)^{+}$ is as in Definition \ref{def:sigmapmais}.

Using the description of the kernel of $\gamma$ of Proposition \ref{pro:kernelgamma},
it is easy to conclude that 
\[\mbox{Im}(f)=\left\{ \varphi\in C^{\infty}\left(\mathbb{S}_{++}^{n}\times\mathbb{S}^{n-1}\cup\mathbb{S}_{+}^{n-1}\times\mathbb{S}_{+}^{n}\right)_{tr},\,\varphi(z',0,z_{n+1},w)=0\right\}.\]

Now we note that if \[
\left(\begin{array}{cc}
P_{+}+G & K\\
T & S\end{array}\right)\in\mathcal{K},\]
then, by Corollary \ref{cor:teo12}, $\sigma(p)^{+}\equiv0$. Hence we can define
the injective function $\tilde{f}:\frac{\mbox{ker}(\gamma)+\mathcal{K}}{\mathcal{K}}\to C^{\infty}\left(\mathbb{S}_{++}^{n}\times\mathbb{S}^{n-1}\cup\mathbb{S}_{+}^{n-1}\times\mathbb{S}_{+}^{n}\right)_{tr}$
by\[
\tilde{f}\left(\left(\begin{array}{cc}
P_{+}+G & K\\
T & S\end{array}\right)+\mathcal{K}\right)=\sigma(p)^{+}.\]

However, again by Corollary \ref{cor:teo12},
\[
\sup_{(z,w)\in \mathbb{S}_{++}^{n}\times\mathbb{S}^{n-1}\cup\mathbb{S}_{+}^{n-1}\times\mathbb{S}_{+}^{n}}\left|\tilde{f}\left(\left(\begin{array}{cc}
P_{+}+G & K\\
T & S\end{array}\right)+\mathcal{K}\right)(z,w)\right|=\]
\[
\left\Vert\sigma(p)^{+}\right\Vert \le\inf_{C\in\mathcal{K}(L^{2}(\mathbb{R}_{+}^{n})\oplus L^{2}(\mathbb{R}^{n-1}))}\left\Vert \left(\begin{array}{cc}
P_{+}+G & K\\
T & S\end{array}\right)+C\right\Vert _{\mathcal{B}(L^{2}(\mathbb{R}_{+}^{n})\oplus L^{2}(\mathbb{R}^{n-1}))}.\]

Therefore $\tilde{f}$ can be extended uniquely to a continuous function \[
j:\frac{\mbox{ker}(\overline{\gamma})}{\mathcal{K}}\to C(\mathbb{S}_{++}^{n}\times\mathbb{S}^{n-1}\cup\mathbb{S}_{+}^{n-1}\times\mathbb{S}_{+}^{n}),\]
as $\mbox{ker}(\overline{\gamma})$ is equal to the closure of $\mbox{ker}(\gamma)$
and contains the compact operators. Clearly $j(A+\mathcal{K})=f(A)$
if $A\in\mbox{ker}(\gamma)$. We have finally to ask ourselves about
the image of $j$. We know that $\mbox{Im}\left(j\right)=\overline{\mbox{Im}\left(f\right)}$, as $j$ is a homomorphism of $C^{*}$-algebras. We know that $\mbox{Im}\left(f\right)$ contains the functions $C^{\infty}(\mathbb{S}_{++}^{n}\times\mathbb{S}^{n-1}\cup\mathbb{S}_{+}^{n-1}\times\mathbb{S}_{+}^{n})$
that are zero in a neighborhood of $z_{n}=0$, because these functions
correspond to principal symbols that are zero when $x$ belongs to a neighborhood
of $\partial\left(\overline{\mathbb{R}^{n}_{+}}\right)$ and, therefore, they satisfy the transmission property. This set
is also contained in the set of functions $C^{\infty}(\mathbb{S}_{++}^{n}\times\mathbb{S}^{n-1}\cup\mathbb{S}_{+}^{n-1}\times\mathbb{S}_{+}^{n})$
that are zero when $z_{n}=0$. The closure of the both sets are equal.
Hence we conclude that the image of $j$ is the set of functions $\varphi\in C(\mathbb{S}_{++}^{n}\times\mathbb{S}^{n-1}\cup\mathbb{S}_{+}^{n-1}\times\mathbb{S}_{+}^{n})$
that are zero for $z_{n}=0$.
\end{proof}
Our second step will be the study of the $K$-theory of the $C^{*}$-algebra $\mbox{Im}(\overline{\gamma})$, which is isomorphic to $\frac{\mathfrak{A}}{\mathcal{J}}$.
Actually it will be enough to obtain certain isomorphisms of the $K$-groups of this algebra.

We start defining a sub algebra of the Wiener-Hopf algebra (see
for instance \cite{RempelSchulze}).
\begin{defn}
Let us denote by $\mathfrak{M}_{0}$ the algebra of bounded
operators on $L^{2}(\mathbb{R}_{+})\oplus\mathbb{C}$ given by 
matrices of the form\[
\left(\begin{array}{cc}
A_{11} & A_{12}\\
A_{21} & A_{22}\end{array}\right),\]
where $A_{11}\in\mathfrak{T}_{0}$, $A_{12}\in L^{2}(\mathbb{R}_{+})$,
which means $A_{12}:\mathbb{C}\to L^{2}(\mathbb{R}_{+})$ is given
by the multiplication of a function on $L^{2}(\mathbb{R}_{+})$, $A_{21}\in L^{2}(\mathbb{R}_{+})^{*}$,
which means that $A_{21}:L^{2}(\mathbb{R}_{+})\to\mathbb{C}$ is a
linear functional on $L^{2}(\mathbb{R}_{+})$, and $A_{22}\in\mathbb{C}$.

In Section \ref{sec:image} we proved that 
\begin{equation}
\label{eq:imgamma}
\mbox{Im}(\overline{\gamma})=C\left(\partial\left(\mathbb{S}^{n-1}_{+}\times\mathbb{S}^{n-1}_{+}\right),\mathfrak{M}_{0})\oplus C(\mathbb{S}_{+x'}^{n-1}\right),
\end{equation}
where we identify $C(\mathbb{S}_{+x'}^{n-1})$ with \[
C(\mathbb{S}_{+x'}^{n-1})\otimes\left(\begin{array}{cc}
\mathbb{C} & 0\\
0 & 0\end{array}\right).\]

\end{defn}
It is well known (see \cite[Lemma 7]{MeloSchrohe} for instance)
that $K_{0}(\mathfrak{M}_{0})=K_{1}(\mathfrak{M}_{0})=0$.
Hence using that $C\left(\partial\left(\mathbb{S}^{n-1}_{+}\times\mathbb{S}^{n-1}_{+}\right),\mathfrak{M}_{0}\right)=C\left(\partial\left(\mathbb{S}^{n-1}_{+}\times\mathbb{S}^{n-1}_{+}\right)\right)\hat{\otimes}\mathfrak{M}_{0}$
and K\"unneth formula, we conclude that the $K$-groups of $C\left(\partial\left(\mathbb{S}^{n-1}_{+}\times\mathbb{S}^{n-1}_{+}\right),\mathfrak{M}_{0}\right)$
are zero.

Let us now define a map $b:C^{\infty}(\mathbb{S}_{+x'}^{n-1})\to\mbox{Im}(\gamma)$
which will induce a $K$-group isomorphism.
\begin{defn}
We define $b:C(\mathbb{S}_{+x'}^{n-1})\to\mbox{Im}(\overline{\gamma})$
by \[
b(g)(z,w)=g(z)I_{L^{2}(\mathbb{R}_{+})\oplus\mathbb{C}},\ (z,w)\in\partial\left(\mathbb{S}^{n-1}_{+}\times\mathbb{S}^{n-1}_{+}\right),\]
where $I_{L^{2}(\mathbb{R}_{+})\oplus\mathbb{C}}$ is the identity
of $\mathcal{B}(L^{2}(\mathbb{R}_{+})\oplus\mathbb{C})$.
\end{defn}

This definition can be made more explicit when $g\in C^{\infty}(\mathbb{S}_{+}^{n-1})$. We choose a function $\chi\in C^{\infty}_{c}(\mathbb{R})$ that is 1 on a neighborhood of 0. Then $b(g)$ is the boundary principal symbol of the multiplication operator by the following function:

 \[\left(\begin{array}{cc}
g\left(\frac{x'}{\left\langle x'\right\rangle },\frac{1}{\left\langle x'\right\rangle }\right)\chi\left(\frac{x_{n}}{\left\langle x'\right\rangle }\right) & 0\\
0 & g\left(\frac{x'}{\left\langle x'\right\rangle },\frac{1}{\left\langle x'\right\rangle }\right)\end{array}\right).\] This definition clearly does not depend on the function $\chi$.

\begin{prop}
\label{pro:biso}
(Analogous to \cite[Corollary 8]{MeloSchrohe}) The $*$-homomorphism $b$ induces isomorphisms $K_{i}(b):K_{i}(C(\mathbb{S}_{+}^{n-1}))\to K_{i}(\mbox{Im}(\overline{\gamma}))$, $i=0,1$.\end{prop}
\begin{proof}
If $F\in C\left(\partial\left(\mathbb{S}^{n-1}_{+}\times\mathbb{S}^{n-1}_{+}\right),\mathfrak{M_{0}}\right)\cap \mbox{Im}(b)$, then we see that the entry 1-1 of the
matrix $F$ can be written as $(z,w)\mapsto g(z)I_{L^{2}(\mathbb{R}_{+}^{n})}$,
where $g\in C(\mathbb{S}_{+x'}^{n-1})$ and $\left( (z,w)\mapsto g(z)I_{L^{2}(\mathbb{R}_{+}^{n})}\right) \in C\left(\partial\left(\mathbb{S}^{n-1}_{+}\times\mathbb{S}^{n-1}_{+}\right),\mathfrak{T_{0}}\right)$.
This implies that $g$ must be zero. Hence, by equation (\ref{eq:imgamma}), we conclude that
\[
\begin{array}{c}
\mbox{Im}(\overline{\gamma})=C\left(\partial\left(\mathbb{S}^{n-1}_{+}\times\mathbb{S}^{n-1}_{+}\right),\mathfrak{M_{0}}\right)\oplus\mbox{Im}(b)\end{array}.\]

As $C\left(\partial\left(\mathbb{S}^{n-1}_{+}\times\mathbb{S}^{n-1}_{+}\right),\mathfrak{M_{0}}\right)$
is an ideal of $\mbox{Im}(\overline{\gamma})$, we get \[
\mbox{Im}(\overline{\gamma})/C\left(\partial\left(\mathbb{S}^{n-1}_{+}\times\mathbb{S}^{n-1}_{+}\right),\mathfrak{M_{0}}\right)\cong\mbox{Im}(b).\]

Let $\mathfrak{R}$ denote $C(\partial\left(\mathbb{S}^{n-1}_{+}\times\mathbb{S}^{n-1}_{+}\right),\mathfrak{M_{0}})$.
Using the exact sequence \foreignlanguage{brazil}{\[
0\to\mathfrak{R}\overset{i}{\to}\mbox{Im}(\overline{\gamma})\overset{\pi}{\to}\mbox{Im}(b)\to0,\]
}where $i$ is the inclusion and $\pi$ is the canonical projection, we obtain the standard six term exact sequence \foreignlanguage{brazil}{\[
\begin{array}{ccccc}
K_{0}(\mathfrak{R}) & \overset{K_{0}(i)}{\to} & K_{0}(\mbox{Im}(\overline{\gamma})) & \overset{K_{0}(\pi)}{\to} & K_{0}(\mbox{Im}(b))\\
\uparrow &  &  &  & \downarrow\\
K_{1}(\mbox{Im}(b)) & \overset{K_{1}(\pi)}{\leftarrow} & K_{1}(\mbox{Im}(\overline{\gamma})) & \overset{K_{1}(i)}{\leftarrow} & K_{1}(\mathfrak{R})\end{array}.\]
}

As $K_{0}(\mathfrak{R})=K_{1}(\mathfrak{R})=0$, we obtain that $K_{i}(\pi):K_{i}(\mbox{Im}(\overline{\gamma}))\to K_{i}(\mbox{Im}(b))$
is an isomorphism, for $i=0,1$. As $b$ is injective, we conclude that $\pi\circ b:C(\mathbb{S}_{+}^{n-1})\to\mbox{Im}(\overline{\gamma})\to\mbox{Im}(b)$
is an isomorphism of $C^{*}$-algebras. Hence $K_{i}(\pi\circ b):K_{i}(C(\mathbb{S}_{+}^{n-1}))\to K_{i}(\mbox{Im}(b))$
is also an isomorphism. We conclude that\[
K_{i}(\pi)^{-1}K_{i}(\pi\circ b):K_{i}(C(\mathbb{S}_{+}^{n-1}))\to K_{i}(\mbox{Im}(\overline{\gamma}))\]
 is an isomorphism. Hence \[
K_{i}(b):K_{i}(C(\mathbb{S}_{+}^{n-1}))\to K_{i}(\mbox{Im}(\overline{\gamma}))\]
is an isomorphism.
\end{proof}
We are now ready to conclude our computation of the $K$-groups of $\mathfrak{A}/\mathcal{K}$. We need some definitions:
\begin{defn}
1) $C_{0\, on\, x_{n}=0}(\mathbb{S}_{++}^{n})$ is the set of continuous
functions on $\mathbb{S}_{++}^{n}\subset\mathbb{R}_{z}^{n+1}$ that
are zero when $z_{n}=0$.

2) $C_{0}(\mathbb{S}_{++}^{n}\times\mathbb{S}_{+}^{n})$ is the
set of continuous functions on $\mathbb{S}_{++z}^{n}\times\mathbb{S}_{+w}^{n}$
that are zero when $z_{n+1}=0$ or $w_{n+1}=0$ or $z_{n}=0$, that is, functions that vanish on the boundary points.

3) $C_{0\, on\, x_{n}=0}(\mathbb{S}_{++}^{n}\times\mathbb{S}_{+}^{n})$
is the set of continuous functions on $\mathbb{S}_{++z}^{n}\times\mathbb{S}_{+w}^{n}$
that are zero when $z_{n}=0$.

4) $C_{0\, on\, x_{n}=0}(\mathbb{S}_{++}^{n}\times\mathbb{S}^{n-1}\cup\mathbb{S}_{+}^{n-1}\times\mathbb{S}_{+}^{n})$
is the set of continuous functions on $\mathbb{S}_{++z}^{n}\times\mathbb{S}_{w}^{n-1}\cup\mathbb{S}_{+z}^{n-1}\times\mathbb{S}_{+w}^{n}$
that are zero when $z_{n}=0$.
\end{defn}
We can use these classes of functions to define the following exact
sequence.\foreignlanguage{brazil}{
\begin{equation}
0\to C_{0}(\mathbb{S}_{++}^{n}\times\mathbb{S}_{+}^{n})\overset{i}{\to}C_{0\, on\, x_{n}=0}(\mathbb{S}_{++}^{n}\times\mathbb{S}_{+}^{n})\overset{\pi}{\to}C_{0\, on\, x_{n}=0}(\mathbb{S}_{++}^{n}\times\mathbb{S}^{n-1}\cup\mathbb{S}_{+}^{n-1}\times\mathbb{S}_{+}^{n})\to0.
\label{eq:ipi}
\end{equation}
}where $i$ is the inclusion and $\pi$ is the restriction of the
functions to the set where $z_{n+1}=0$ or $w_{n+1}=0$. We also need
to use the following maps.
\begin{defn}
1) The map $m':C_{0\, on\, x_{n}=0}(\mathbb{S}_{++}^{n})\to C_{0\, on\, x_{n}=0}(\mathbb{S}_{++}^{n}\times\mathbb{S}^{n-1}\cup\mathbb{S}_{+}^{n-1}\times\mathbb{S}_{+}^{n})$
is defined by\[
m'(f)(z,w)=f(z).\]

2) The map $m'':C_{0\, on\, x_{n}=0}(\mathbb{S}_{++}^{n})\to C_{0\, on\, x_{n}=0}(\mathbb{S}_{++}^{n}\times\mathbb{S}_{+}^{n})$
is defined by\[
m''(f)(z,w)=f(z)\]
Notice that $m'(f)$ is a restriction of $m''(f)$.

3) Let us choose $p\in\mathbb{S}^{n-1}$, then we can define $s:C_{0\, on\, x_{n}=0}(\mathbb{S}_{++}^{n}\times\mathbb{S}^{n-1}\cup\mathbb{S}_{+}^{n-1}\times\mathbb{S}_{+}^{n})\to C_{0\, on\, x_{n}=0}(\mathbb{S}_{++}^{n})$
by \[
s(f)(z)=f(z,p).\]

It is clear that $s\circ m'$ is equal to the identity.\end{defn}
\begin{prop}
\label{pro:hi}
(Analogous to  \cite[Proposition 11]{MeloSchrohe}). For $i=0,1$,
there is an isomorphism
\[
h_{i}:K_{i}(C_{0\, on\, x_{n}=0}(\mathbb{S}_{++}^{n}\times\mathbb{S}^{n-1}\cup\mathbb{S}_{+}^{n-1}\times\mathbb{S}_{+}^{n}))\to K_{i}(C_{0\, on\, x_{n}=0}(\mathbb{S}_{++}^{n}))\oplus K_{1-i}(C_{0}(\mathbb{S}_{++}^{n}\times\mathbb{S}_{+}^{n})).\]
This isomorphism $h_{i}$ is such that the injection \[
I_{i}:K_{i}(C_{0\, on\, x_{n}=0}(\mathbb{S}_{++}^{n}))\to K_{i}(C_{0\, on\, x_{n}=0}(\mathbb{S}_{++}^{n}))\oplus K_{1-i}(C_{0}(\mathbb{S}_{++}^{n}\times\mathbb{S}_{+}^{n}))\]
 is given by $h_{i}\circ K_{i}(m')$, where $K_{i}(m')$ is the group homomorphism induced by $m'$. In particular $K_{i}(m')$ is injective.
Moreover the projection \[
p:K_{i}(C_{0\, on\, x_{n}=0}(\mathbb{S}_{++}^{n}))\oplus K_{1-i}(C_{0}(\mathbb{S}_{++}^{n}\times\mathbb{S}_{+}^{n}))\to K_{1-i}(C_{0}(\mathbb{S}_{++}^{n}\times\mathbb{S}_{+}^{n}))\]
 is given by $h_{i}\circ\delta$, where $\delta$ is the connecting map
associated to the six term exact sequence provided by the exact sequence (\ref{eq:ipi}).
\end{prop}
\begin{proof}
As $\mathbb{S}_{+}^{n}$ can be deformed continuously to a point,
we conclude that the function $m''$ induces isomorphisms, for $i=0,1$,\[
K_{i}(m''):K_{i}\left(C_{0\, on\, x_{n}=0}(\mathbb{S}_{++}^{n})\right)\to K_{i}\left(C_{0\, on\, x_{n}=0}(\mathbb{S}_{++}^{n}\times\mathbb{S}_{+}^{n})\right).\]

We see easily that $\pi\circ m''=m'$. Hence, denoting $C_{0\, on\, x_{n}=0}(\mathbb{S}_{++}^{n}\times\mathbb{S}^{n-1}\cup\mathbb{S}_{+}^{n-1}\times\mathbb{S}_{+}^{n})$
just by $C_{0\, on\, x_{n}}^{\partial}$, we can rewrite the standart
six term exact sequence associated to (\ref{eq:ipi}) \foreignlanguage{brazil}{\[
\begin{array}{ccccc}
K_{0}\left(C_{0}(\mathbb{S}_{++}^{n}\times\mathbb{S}_{+}^{n})\right) & \overset{K_{0}(m'')^{-1}\circ K_{0}(i)}{\longrightarrow} & K_{0}\left(C_{0\, on\, x_{n}=0}(\mathbb{S}_{++}^{n})\right) & \overset{K_{0}(m')}{\longrightarrow} & K_{0}\left(C_{0\, on\, x_{n}}^{\partial}\right)\\
\uparrow &  &  &  & \downarrow\\
K_{1}\left(C_{0\, on\, x_{n}}^{\partial}\right) & \overset{K_{1}(m')}{\longleftarrow} & K_{1}\left(C_{0\, on\, x_{n}=0}(\mathbb{S}_{++}^{n})\right) & \overset{K_{1}(m'')^{-1}\circ K_{1}(i)}{\longleftarrow} & K_{1}\left(C_{0}(\mathbb{S}_{++}^{n}\times\mathbb{S}_{+}^{n})\right)\end{array}.\]
}

However $s\circ m'$ is the identity. Therefore $K_{i}(s)K_{i}(m')=K_{i}(id)=id_{K_{i}}$
and $K_{i}(m')$ is injective, which implies that $\mbox{ker}K_{i}(m')=\{0\}=\mbox{Im}(K_{i}(m'')^{-1}\circ K_{i}(i))$. We conclude that the
following sequence are split exact: 
\[\begin{array}{c}
0\to K_{i}\left(C_{0\, on\, x_{n}=0}(S_{++}^{n})\right)\begin{array}{c}
\overset{K_{i}(m')}{\to}\\
\underset{K_{i}(s)}{\leftarrow}\end{array}K_{i}\left(C_{0\, on\, x_{n}=0}^{\partial}\right)\overset{\delta}{\longrightarrow}K_{1-i}\left(C_{0}(S_{++}^{n}\times S_{+}^{n})\right)\rightarrow 0.\end{array}\]

Now the result follows easily.\end{proof}
\begin{thm}
(Analogous of \cite[Theorem 3]{MeloSchrohe}) Let $r:C(\mathbb{S}_{++}^{n})\to C(\mathbb{S}_{+}^{n-1})$
be the restriction to $z_{n}=0$, where $z\in\mathbb{S}^{n}_{++}\subset\mathbb{R}^{n+1}$ and let, for $i=0,1$, 
$K_{i}(r):K_{i}\left(C(\mathbb{S}_{++}^{n})\right)\to K_{i}\left(C(\mathbb{S}_{+}^{n-1})\right)$ denote the induced homomorphisms. There is an exact sequence\[
0\to ker(K_{i}(r))\oplus K_{1-i}(C_{0}(\mathbb{S}_{++}^{n}\times\mathbb{S}_{+}^{n}))\to K_{i}\left(\frac{\mathfrak{A}}{\mathcal{K}}\right)\to Im(K_{i}(r))\to0.\]
\end{thm}
\begin{proof}
Continuous functions in $C(\mathbb{S}_{++}^{n})$ can be canonically viewed as multiplication operators. By taking the class, modulo $\mathcal{K}$, of such an operator, one defines a homomorphism $m:C(\mathbb{S}_{++}^{n})\to\mathfrak{A}/\mathcal{K}$. This can be made more precise on the dense subalgebra $C^{\infty}(\mathbb{S}_{++}^{n})$: For each $f\in C^{\infty}(\mathbb{S}_{++}^{n})$, we assign a function
 $g\in S_{cl}^{0,0}(\mathbb{R}^{n}\times \mathbb{R}^{n})$ given by $g(x,\xi)=f\circ RC(x)$, for $x_{n}>0$. Then \[
m(f)=\left(\begin{array}{cc}
op(g)_{+} & 0\\
0 & g|_{\mathbb{R}_{x'}^{n-1}\times \mathbb{R}_{\xi'}^{n-1}}\end{array}\right)+\mathcal{K}.\]
This clearly does not depend on the values of $g$ for $x_{n}<0$ and, therefore, it is well defined.

Now let us consider the following commutative diagram of $C^{*}$-algebras homomorphisms:
\begin{equation}
\begin{array}{ccccccccc}
0 & \to & \mathcal{J}/\mathcal{K} & \to & \mathfrak{A}/\mathcal{K} & \to & \mathfrak{A}/\mathcal{J} & \to & 0\\
 &  & \uparrow m &  & \uparrow m &  & \uparrow b\\
0 & \to & C_{0\, on\, x_{n}=0}(\mathbb{S}_{++}^{n}) & \overset{i}{\longrightarrow} & C(\mathbb{S}_{++}^{n}) & \overset{r}{\longrightarrow} & C(\mathbb{S}_{+}^{n-1}) & \to & 0\end{array},
\label{eq:seq2}
\end{equation}
where the upper row is the sequence in (\ref{eq:basica}) and the lower row is the exact sequence defined by the inclusion $i$ and the restriction map to $z_{n}=0$, denoted by $r$.

Let us denote by $\delta$ and $\exp$ the index and the exponential
maps associated to the sequence (\ref{eq:basica}). By $\delta^{0}$ and $\exp^{0}$ we denote the index and the
exponential maps associated to the exact sequence of the lower row of the diagram (\ref{eq:seq2}).

As these maps are natural, we obtain two commutative diagrams\foreignlanguage{brazil}{\[
\begin{array}{ccc}
K_{1}\left(\frac{\mathfrak{A}}{\mathcal{J}}\right) & \overset{\delta}{\longrightarrow} & K_{0}\left(\frac{\mathcal{J}}{\mathcal{K}}\right)\\
\uparrow K_{1}(b) &  & \uparrow K_{0}(m)\\
K_{1}(C(\mathbb{S}_{+}^{n-1})) & \overset{\delta^{0}}{\longrightarrow} & K_{0}(C_{0\, on\, x_{n}=0}(\mathbb{S}_{++}^{n}))\end{array}\mbox{and}\begin{array}{ccc}
K_{0}\left(\frac{\mathfrak{A}}{\mathcal{J}}\right) & \overset{\exp}{\longrightarrow} & K_{1}\left(\frac{\mathcal{J}}{\mathcal{K}}\right)\\
\uparrow K_{0}(b) &  & \uparrow K_{1}(m)\\
K_{0}(C(\mathbb{S}_{+}^{n-1})) & \overset{\exp^{0}}{\longrightarrow} & K_{1}(C_{0\, on\, x_{n}=0}(\mathbb{S}_{++}^{n}))\end{array}\]
}

Let $j:\mathcal{J}/\mathcal{K}\to C_{0\, on\, x_{n}=0}(\mathbb{S}_{++}^{n}\times\mathbb{S}^{n-1}\cup\mathbb{S}_{+}^{n-1}\times\mathbb{S}_{+}^{n})$
be the isomorphism given by Proposition \foreignlanguage{brazil}{\ref{pro:kergamma em simb prin de p}}.
It is clear that $m'=\left.j\circ m\right|_{C_{0\, on\, x_{n}=0}(\mathbb{S}_{++}^{n})}$.
We saw in the proof of Proposition \ref{pro:hi} that $K_{i}(m')$ is injective. As $K_{i}(j)$ is an isomorphism, we conclude that $K_{i}(m):K_{i}(C_{0\, on\, x_{n}=0}(\mathbb{S}_{++}^{n}))\to K_{i}\left(\frac{\mathcal{J}}{\mathcal{K}}\right)$
is also injective. Therefore, as $K_{i}(b)$ is an isomorphism (Proposition \ref{pro:biso}), we conclude that $K_{1}(b)^{-1}\left(\ker(\delta)\right)=\ker(\delta^{0})$
and $K_{0}(b)^{-1}\left(\ker(\exp)\right)=\ker(\exp^{0})$.

Using the isomorphisms $K_{i}(b)$ and $K_{i}(j)$, we define $\mbox{exp}''$ and $\delta ''$ by \foreignlanguage{brazil}{$\delta''=K_{0}(j)\circ\delta\circ K_{1}(b)=K_{0}(j)\circ K_{0}(m)\circ\delta^{0}$}
and \foreignlanguage{brazil}{$\exp''=K_{1}(j)\circ\exp\circ K_{0}(b)=K_{1}(j)\circ K_{1}(m)\circ\exp^{0}$}.
We note that \foreignlanguage{brazil}{$\ker(\delta'')=K_{1}(b)^{-1}\left(\ker(\delta)\right)=\ker(\delta^{0})$}
and \foreignlanguage{brazil}{$\ker(\exp'')=K_{0}(b)^{-1}\left(\ker(\exp)\right)=\ker(\exp^{0})$}.
Hence the six term exact sequence associated to (\ref{eq:basica}) can be rewritten as:  \[
\begin{array}{ccccc}
K_{0}\left(C_{0\, on\, x_{n}=0}(\mathbb{S}_{++}^{n}\times\mathbb{S}^{n-1}\cup\mathbb{S}_{+}^{n-1}\times\mathbb{S}_{+}^{n})\right) & \to & K_{0}\left(\frac{\mathfrak{A}}{\mathcal{K}}\right) & \to & K_{0}\left(C(\mathbb{S}_{+}^{n-1})\right)\\
\uparrow\delta'' &  &  &  & \downarrow\exp''\\
K_{1}\left(C(\mathbb{S}_{+}^{n-1})\right) & \leftarrow & K_{1}\left(\frac{\mathfrak{A}}{\mathcal{K}}\right) & \leftarrow & K_{1}\left(C_{0\, on\, x_{n}=0}(\mathbb{S}_{++}^{n}\times\mathbb{S}^{n-1}\cup\mathbb{S}_{+}^{n-1}\times\mathbb{S}_{+}^{n})\right)\end{array}.\]

Taking quotients, we obtain from this exact sequence that
\begin{equation}
0\to\frac{K_{0}(C_{0\, on\, x_{n}=0}(\mathbb{S}_{++}^{n}\times\mathbb{S}^{n-1}\cup\mathbb{S}_{+}^{n-1}\times\mathbb{S}_{+}^{n}))}{\mbox{Im}(\delta'')}\to K_{0}\left(\frac{\mathfrak{A}}{\mathcal{K}}\right)\to\ker(\exp'')\to0,
\label{eq:A}
\end{equation}

\begin{equation}
0\to\frac{K_{1}(C_{0\, on\, x_{n}=0}(\mathbb{S}_{++}^{n}\times\mathbb{S}^{n-1}\cup\mathbb{S}_{+}^{n-1}\times\mathbb{S}_{+}^{n}))}{\mbox{Im}(\exp'')}\to K_{1}\left(\frac{\mathfrak{A}}{\mathcal{K}}\right)\to\ker(\delta'')\to0.
\label{eq:B}
\end{equation}

Using the standart exact sequence associated to (\ref{eq:seq2}), we obtain the following isomorphisms:\[
\begin{array}{c}
\frac{K_{0}(C_{0\, on\, x_{n}=0}(\mathbb{S}_{++}^{n}))}{\mbox{Im}(\delta^{0})}=\frac{K_{0}(C_{0\, on\, x_{n}=0}(\mathbb{S}_{++}^{n}))}{\mbox{ker}(K_{0}(i))}\cong\mbox{Im}(K_{0}(i))=\ker K_{o}(r),\\
\frac{K_{1}(C_{0\, on\, x_{n}=0}(\mathbb{S}_{++}^{n}))}{\mbox{Im}(\exp^{0})}=\frac{K_{1}(C_{0\, on\, x_{n}=0}(\mathbb{S}_{++}^{n}))}{\mbox{ker}(K_{1}(i))}\cong\mbox{Im}(K_{1}(i))=\ker K_{1}(r).\end{array}\]

As $\delta''=K_{0}(j)\circ\delta\circ K_{1}(b)=K_{0}(m')\circ\delta^{0}$
and $\exp''=K_{1}(j)\circ K_{1}(m)\circ\exp^{0}=K_{1}(m')\circ\exp^{0}$,
we conclude that $\mbox{Im}(\delta'')=\mbox{Im}(K_{0}(m')\circ\delta^{0})$
and $\mbox{Im}(\exp'')=\mbox{Im}(K_{1}(m')\circ\exp^{0})$. Using the maps
 $I_{i}$ and $h_{i}$ defined in Proposition \ref{pro:hi},
we obtain:\[
\begin{array}{c}
\frac{K_{0}(C_{0\, on\, x_{n}=0}(\mathbb{S}_{++}^{n}\times\mathbb{S}^{n-1}\cup\mathbb{S}_{+}^{n-1}\times\mathbb{S}_{+}^{n}))}{\mbox{Im}(\delta'')}\cong\frac{h_{0}\left( K_{0}(C_{0\, on\, x_{n}=0}(\mathbb{S}_{++}^{n}\times\mathbb{S}^{n-1}\cup\mathbb{S}_{+}^{n-1}\times\mathbb{S}_{+}^{n}))\right)}{\mbox{Im}(h_{0}\circ K_{0}(m')\circ\delta^{0})}\cong\frac{K_{0}(C_{0\, on\, x_{n}=0}(\mathbb{S}_{++}^{n}))\oplus K_{1}(C_{0}(\mathbb{S}_{++}^{n}\times\mathbb{S}_{+}^{n}))}{\mbox{Im}(I_{0}\circ\delta^{0})}\\
 \cong K_{1}(C_{0}(\mathbb{S}_{++}^{n}\times\mathbb{S}_{+}^{n}))\oplus\frac{K_{0}(C_{0\, on\, x_{n}=0}(\mathbb{S}_{++}^{n}))}{\mbox{Im}(\delta^{0})}\cong K_{1}(C_{0}(\mathbb{S}_{++}^{n}\times\mathbb{S}_{+}^{n}))\oplus\ker K_{0}(r).\end{array}\]

and\[
\begin{array}{c}
\frac{K_{1}(C_{0\, on\, x_{n}=0}(\mathbb{S}_{++}^{n}\times\mathbb{S}^{n-1}\cup\mathbb{S}_{+}^{n-1}\times\mathbb{S}_{+}^{n}))}{\mbox{Im}(\exp'')}\cong\frac{h_{1}\left(K_{1}(C_{0\, on\, x_{n}=0}(\mathbb{S}_{++}^{n}\times\mathbb{S}^{n-1}\cup\mathbb{S}_{+}^{n-1}\times\mathbb{S}_{+}^{n}))\right)}{\mbox{Im}(h_{1}\circ K_{1}(m')\circ\exp^{0})}\cong 
\frac{K_{1}(C_{0\, on\, x_{n}=0}(\mathbb{S}_{++}^{n}))\oplus K_{0}(C_{0}(\mathbb{S}_{++}^{n}\times\mathbb{S}_{+}^{n}))}{\mbox{Im}(I_{1}\circ\exp^{0})}\\ \cong K_{0}(C_{0}(\mathbb{S}_{++}^{n}\times\mathbb{S}_{+}^{n}))\oplus\frac{K_{1}(C_{0\, on\, x_{n}=0}(\mathbb{S}_{++}^{n}))}{\mbox{Im}(\exp^{0})}\cong K_{0}(C_{0}(\mathbb{S}_{++}^{n}\times\mathbb{S}_{+}^{n}))\oplus\ker K_{1}(r).\end{array}\]

The exact sequence in (\ref{eq:A})
turns into \[
0\to K_{1}(C_{0}(\mathbb{S}_{++}^{n}\times\mathbb{S}_{+}^{n}))\oplus\ker K_{0}(r)\to K_{0}\left(\frac{\mathfrak{A}}{\mathcal{K}}\right)\to\ker\exp''\to0,\]
as $\ker(\exp'')=\ker(\exp^{0})\cong\mbox{Im}K_{0}(r)$, we finally
obtain the following exact sequence:\[
0\to K_{1}(C_{0}(\mathbb{S}_{++}^{n}\times\mathbb{S}_{+}^{n}))\oplus\ker K_{0}(r)\to K_{0}\left(\frac{\mathfrak{A}}{\mathcal{K}}\right)\to\mbox{Im}K_{0}(r)\to0.\]
Similarly, using sequence (\ref{eq:B}) and $\ker(\delta'')=\ker(\delta^{0})\cong\mbox{Im}K_{1}(r)$,
we conclude that \[
0\to K_{0}(C_{0}(\mathbb{S}_{++}^{n}\times\mathbb{S}_{+}^{n}))\oplus\ker K_{1}(r)\to K_{1}\left(\frac{\mathfrak{A}}{\mathcal{K}}\right)\to\mbox{Im}K_{1}(r)\to0.\]
\end{proof}
\begin{cor}
The $K$-groups of $\mathfrak{A}/\mathcal{K}$ are given by \[
K_{i}(\mathfrak{A}/\mathcal{K})\cong K_{i}(C(\mathbb{S}_{++}^{n}))\oplus K_{1-i}(C_{0}(\mathbb{S}_{++}^{n}\times\mathbb{S}_{+}^{n})).\]
\end{cor}
\begin{proof}
As $K_{i}\left(C(\mathbb{S}_{+}^{n-1})\right)$ is a free abelian group,
once $\mathbb{S}_{+}^{n-1}$ is a contractible space, see \cite{Rordamkteoria}, we conclude that $\mbox{Im}(K_{i}(r))$ is an free abelian
group.
Hence both sequences below split: \[
0\to ker(K_{i}(r))\oplus K_{1-i}(C_{0}(\mathbb{S}_{++}^{n}\times\mathbb{S}_{+}^{n}))\to K_{i}\left(\frac{\mathfrak{A}}{\mathcal{K}}\right)\to Im(K_{i}(r))\to0\]
and

\[
0\to ker(K_{i}(r))\to K_{i}\left(C(\mathbb{S}_{++}^{n})\right)\to Im(K_{i}(r))\to0.\]
Therefore\[
K_{i}\left(\frac{\mathfrak{A}}{\mathcal{K}}\right)\cong Im(K_{i}(r))\oplus ker(K_{i}(r))\oplus K_{1-i}(C_{0}(\mathbb{S}_{++}^{n}\times\mathbb{S}_{+}^{n}))\]
\[
K_{i}\left(C(\mathbb{S}_{++}^{n})\right)\cong Im(K_{i}(r))\oplus ker(K_{i}(r)).\]

\end{proof}
Finally we obtain our main result.
\begin{thm}
The $K$-groups of $\mathfrak{A}/\mathcal{K}$ are\[
\begin{array}{c}
K_{0}\left(\mathfrak{A}/\mathcal{K}\right)=\mathbb{Z}\\
K_{1}\left(\mathfrak{A}/\mathcal{K}\right)=\mathbb{Z}.\end{array}\]
\end{thm}
\begin{proof}
As $S_{+}^{n-1}$ is a contractible set, we conclude by \cite{Rordamkteoria}
that\[
\begin{array}{c}
K_{0}\left(C(\mathbb{S}_{+}^{n-1})\right)=\mathbb{Z}\\
K_{1}\left(C(\mathbb{S}_{+}^{n-1})\right)=0.\end{array}\]

Regarding $C_{0}(\mathbb{S}_{++}^{n}\times\mathbb{S}_{+}^{n})$,
we observe that $\mathbb{S}_{++}^{n}$ and $\mathbb{S}_{+}^{n}$
are both homeomorphic to \[
B_{1}^{n}(0)=\{x\in\mathbb{R}^{n},|x|\le1\}.\]
Therefore $C_{0}(\mathbb{S}_{++}^{n}\times\mathbb{S}_{+}^{n})$
is isomorphic to $C_{0}(B_{1}^{n}(0)\times B_{1}^{n}(0))$, that is,
the continuous functions on $B_{1}^{n}(0)\times B_{1}^{n}(0)$ that
are zero on the boundary points. However using an homomorphism that takes the
interior of $B_{1}^{n}(0)$ to $\mathbb{R}^{n}$, we conclude that $C_{0}(B_{1}^{n}(0)\times B_{1}^{n}(0))$
is isomorphic to $C_{0}(\mathbb{R}^{2n})$ and finally\[
\begin{array}{c}
K_{0}\left(C_{0}(\mathbb{S}_{++}^{n}\times\mathbb{S}_{+}^{n})\right)=K_{0}\left(C_{0}(\mathbb{R}^{2n})\right)=\mathbb{Z}\\
K_{1}\left(C_{0}(\mathbb{S}_{++}^{n}\times\mathbb{S}_{+}^{n})\right)=K_{1}\left(C_{0}(\mathbb{R}^{2n})\right)=0.\end{array}\]

\end{proof}
It is interesting to note that the $K$-groups of our algebra are equal to the $K$-groups of the
algebra of SG-operators of order (0,0) acting in $\mathbb{R}^{n}$,
as it was computed by Nicola \cite{Nicolakteoria}. They are also
the same $K$-groups of the Boutet de Monvel algebra on compact manifolds
of dimension 2, whose genus is equal to $0$ and whose border is connected,
as it was shown in \cite[Section 6]{MeloSchrohe}.

{\section*{Acknowledgements}  The authors would like to thank Elmar Schrohe for many fruitful discussions. Pedro Lopes was supported by the Brazilian agency CNPq, Conselho Nacional de Desenvolvimento Cient\'ifico e Tecnol\'ogico -  (Processo n$^{\circ}$ 142185/2007-8). Severino Melo was also partially supported by CNPq -  (Processo n$^{\circ}$ 304783/2009-9).  

\bibliographystyle{amsplain}

\end{document}